\documentclass[preprint,17pt, sort&compress]{elsarticle}
\usepackage{lineno,hyperref}
\usepackage{upgreek}
\modulolinenumbers[5]
\usepackage[margin=2.5cm]{geometry}
\journal{Journal of \LaTeX\ Templates}
\usepackage{tikz}
\usepackage{circuitikz}

\usepackage{enumitem}

\newtheorem{remark}{\sc \bf Remark}
\newtheorem{definition }{\sc \bf Definition }
\usepackage{mathtools}

\usepackage{algorithm, algorithmic}

\usepackage{amsmath} %
\usepackage{graphicx}
\usepackage{caption}
\usepackage{float}
\usepackage{units}
\usepackage{bm}

\usepackage{graphicx}
\makeatletter
\def\fps@figure{hbtp}
\def\fps@table{hbtp}
\makeatother

\newtheorem{example}{\sc \bf Example}

\newcommand{\proofofref}{}
\newproof{zproofof}{Proof of \proofofref}

\newtheorem{lemma}{\sc \bf Lemma}
\usepackage{units}
\usepackage{bm}
\usepackage{booktabs}
\usepackage{multirow}
\usepackage{amssymb,amsmath}
\usepackage{upgreek}
\usepackage[section]{placeins}
\usepackage{upgreek}
\usepackage{booktabs}
\usepackage{multicol}
\usepackage{multirow}
\usepackage{array}
\usepackage[utf8]{inputenc}

\usepackage{amssymb} %
\usepackage{latexsym} %

\usepackage{graphicx}
\hypersetup{colorlinks,linkcolor={blue},citecolor={blue},urlcolor={red}}
\def\mbx{\mathbf{x}}

\begin{document}

\journal{{\textcolor{blue}{\textbf{\normalfont   }}}}

\begin{frontmatter}
\title{Numerical approximation based on deep convolutional neural network for high-dimensional fully nonlinear merged PDEs and 2BSDEs}

\author[mymainaddress1]{Xu Xiao}
\ead{xiaoxu961004@gmail.com}
\author[mymainaddress1]{Wenlin Qiu\textcolor{blue}{$^\ast$}}
\ead{qwllkx12379@163.com}
\author[mymainaddress2]{Omid Nikan}
\ead{omidnikan77@yahoo.com}

\cortext[mycorrespondingauthor]{Corresponding author.  This work was supported by Postgraduate Scientific Research Innovation Project of Hunan Province (No. CX20220454)}

\address[mymainaddress1]{\normalfont Key Laboratory of Computing and Stochastic Mathematics (Ministry of Education), School of Mathematics and Statistics, Hunan Normal University, Changsha, Hunan 410081, P. R. China}
\address[mymainaddress2]{\normalfont School of Mathematics, Iran University of Science and Technology, Narmak, Tehran 16846-13114, Iran}


\begin{abstract}
This paper proposes two efficient approximation methods to solve high-dimensional fully nonlinear partial differential equations (NPDEs) and second-order backward stochastic differential equations (2BSDEs), where such high-dimensional fully NPDEs are extremely difficult to solve because the computational cost of standard approximation methods grows exponentially with the number of dimensions. Therefore, we consider the following methods to overcome this difficulty. For the merged fully NPDEs and 2BSDEs system, combined with the time forward discretization and ReLU function, we use multi-scale deep learning fusion and convolutional neural network (CNN) techniques to obtain two numerical approximation schemes, respectively. Finally,  three practical high-dimensional test problems involving Allen-Cahn, Black-Scholes-Barentblatt, and Hamiltonian-Jacobi-Bellman equations are given so that the first proposed method exhibits higher efficiency and accuracy than the existing method, while the second proposed method can extend the dimensionality of the completely NPDEs-2BSDEs system over $400$ dimensions, from which the numerical results  highlight the effectiveness of proposed methods.
\end{abstract}

\begin{keyword}
 Convolutional neural network \sep ReLU \sep second-order backward stochastic differential equations \sep high-dimensional problems \sep Allen-Cahn equation \sep Black-Scholes-Barentblatt equation \sep Hamiltonian-Jacobi-Bellman equation \sep numerical experiments
\bigskip
\MSC[2020] 65M22\sep 60H15\sep 65C30\sep 68T07
\end{keyword}
\end{frontmatter}

\section{Introduction}\label{Se1}

\vskip 0.2mm
Nonlinear Partial differential equations (NPDEs)  play a key role in a large number of models, from finance to physics. Objects such as wave functions related to quantum physical systems, value functions which depict the fair prices of financial derivatives in pricing models, or value functions which depict the expected maximum utility in portfolio optimization problems that are usually presented as the solutions of NPDEs.

\vskip 0.2mm
Roughly speaking, the non-linearity in PDEs used in financial engineering above is derived from the trade mix (the trade mix and utility of hedging financial derivatives claims in the case of the derivatives pricing problem must be maximized in the case of the portfolio optimization problem). The authors of \cite{7,45}  adopted derivative pricing models with distinguishing lending rates.  Cr\'{e}pey et al. \cite{24} considered derivative pricing models incorporating the default risk of the issuer of the financial derivative.  The authors of \cite{5} proposed the models for the pricing of financial derivatives on untradable underlyings  and analyzed, e.g., financial derivatives on the temperature or mortality-dependent financial derivatives.
Amadori \cite{1} considered the models incorporating that the trading strategy effects the price processes though the demand and supply.

\vskip 0.2mm
The resulting PDEs from these models are usually high-dimensional, since the associated trading portfolio often involves a whole basket of financial assets (see \cite{7,24}). These high-dimensional NPDEs are often exceedingly difficult to be solved approximately. Furthermore, due to the practical relevance of the aforementioned PDEs, there is a strong demand in the financial engineering industry to approximation solutions to such high-dimensional nonlinear parabolic PDEs.

\vskip 0.2mm
There are lots of numerical approaches for solving parabolic NPDEs approximatively in the literature, from which, some of these methods are deterministic approximations, while others are stochastic approximations that depend on appropriate probabilistic representations of the corresponding PDE solutions, e.g., probabilistic representations in view of backward stochastic differential equations (BSDEs) (see \cite{85,87}), probabilistic representations in view of 2BSDEs (see \cite{22}), probabilistic representations in view of branching diffusions (see \cite{55}), and probabilistic representations in view of extensions of the classical Feynman-Kac formula (see \cite{84}). Then, we can refer to some articles specifically, e.g., deterministic approximation approaches for PDEs (see \cite{69,95}), probabilistic approximation approaches for PDEs based on time discretizations of BSDEs (see \cite{6,7,13,18,19,20,21,26,27,28,32,43,44,45,46,48,56,73,75,81,82,93}), probabilistic approximation approaches for PDEs in view of suitable deep learning approximations for BSDEs (see \cite{33,52}), probabilistic approximation approaches for BSDEs in view of Wiener Chaos expansions (see \cite{14}), probabilistic approximation approaches for BSDEs in view of sparse grid approximations (see \cite{40}), probabilistic approximation approaches for PDEs based on branching diffusion representations (see \cite{17,55}), probabilistic approximation approaches for PDEs in view of time discretization of 2BSDEs (see \cite{12,22,49,62}), etc.

\vskip 0.2mm
However, most of the above approximation techniques are only applicable when the dimension $d$ of PDEs/BSDEs is quite small or only when there are strict constraints on the parameters or the type of PDEs considered (e.g., small nonlinearities, small terminal/initial conditions, the semi-linear structure of PDEs, etc). Therefore, to yield the numerical solutions of high-dimensional nonlinear PDEs, this is still an exceedingly difficult task, and there are only a few cases where practical algorithms for high-dimensional PDEs can be considered (see \cite{29,33,52,55}). Especially, to our knowledge, few practical algorithms for high-dimensional fully nonlinear parabolic PDEs currently exist in the scientific literature.

\vskip 0.2mm
This paper intends to solve this difficulty and present new results, i.e., we solve the fully nonlinear merged PDEs and 2BSDEs with a new algorithm. Regarding the proposed problem, Beck et al. \cite{Beck} first consider that by utilizing some properties from Peng's nonlinear expectation in high-dimensional space (see \cite{91}). The proposed algorithm uses a connection between PDEs and 2BSDEs (see Cheridito et al. \cite{22}) to yield a merged formulation of PDEs and 2BSDEs, whose approximated solutions can be obtained via combining time discretizations with a neural network (NN) based on deep learning (see \cite{8,16,33,52,68,67,69,95}). Loosely speaking, the merged formulation allows us to establish the original partial differential problem as a learning problem. The random loss function for the deep neural network in our method can be given by the error between the prescribed terminal condition of 2BSDEs and the neural network in view of forward time discretization of 2BSDEs. In fact, a corresponding deep-learning approximation algorithm for semilinear-type PDEs in view of forward BSDEs has been recently considered in \cite{33,52}. A crucial distinction between \cite{33,52} and our work is that herein we depend on the connection between fully nonlinear PDEs and 2BSDEs given in \cite{22}, while \cite{33,52} depend on the almost classical combination between PDEs and BSDEs (see \cite{85,87}). Besides, although Beck et al. \cite{Beck} have considered the merged construction of fully nonlinear PDEs and 2BSDEs, there is still room for improvement. Under the limitation of computer memory, since they only consider linear neural networks, they can only calculate general high-dimensional nonlinear parabolic problems and cannot calculate higher-dimensional problems (e.g., more than 200 dimensions), and further the approximated error can also be reduced in terms of computational accuracy. These inspired us to carry out the following research.

\par
 The main contributions of this work are as follows:
(i) we improve the method of Beck et al. \cite{Beck} in order to further improve the accuracy of the solution. We apply multi-scale fusion technology \cite{Hu,Kong,Wang} to the original neural network model, that is, use different scales to spatially discretize it, and finally use the merged results. This paper currently uses 4 scales for fusion, (ii) we also generalize the approach in \cite{Beck} so that higher-dimensional models can be solved. The method of \cite{Beck} is to spatially discretize the time-discrete data in the form of vectors. We first arrange the time-discrete data into a matrix and then use the convolutional neural networks \cite{Long,Ronneberger} for spatial discretizations. From the experimental results, the dimension of the solution is further expanded, and the time spent is also shorter. At present, we mainly enumerate numerical experiments in 256 and 400 dimensions,
(iii) we mainly solve three practical high-dimensional examples, which possess the significant physical background, namely, the Allen-Cahn (AC), the Hamilton-Jacobi-Bellman (HJB), and the Black-Scholes-Barenblatt (BSB) equations. The numerical results can demonstrate the effectiveness of the proposed approximation method and
(iv) the proposed strategy considers advanced optimization algorithms, i.e., Adam optimizer and stochastic gradient descent-type optimization.
\par
Following these ideas,  the organizational structure of this work is as follows.
Section \ref{Se2} introduces merged construction of PDEs and 2BSDEs.
Section  \ref{Se3} presents the forward temporal discretizations of the merged PDEs-2BSDEs system, spatial discretizations based on multiscale deep learning fusion and convolutional neural network, respectively, and corresponding optimization algorithms.
Section  \ref{Se4} reports some experiments for numerical solutions of the merged PDEs-2BSDEs system, concretely, containing the high-dimensional AC, BSB  and HJB equations. Finally, Section  \ref{Se5}  summarizes the concluding remarks.

\section{Merged PDEs-2BSDEs system}\label{Se2}
\vskip 0.2mm
This section mainly intends to obtain a merged PDEs-2BSDEs system. First, we shall introduce the fully nonlinear second-order PDEs. Besides, Table \ref{table:nota} summarizes some notions and notations used in this paper.

\begin{table}[H]
\begin{center} \small
\caption{Summarization of notion and notations.}
\label{table:nota}
\begin{tabular}{|m{1.4cm}<{\centering}|c|c|}
\hline
& Notion & Notation\\
\hline
{Function symbol}
&The needed unkonwn function & $u(t,\mathbf{x})$(abbr. $u$) \\
&The boundary function of time & $\hat{g}(\mathbf{x})(u(T,\mathbf{x})=\hat{g}(\mathbf{x}))$\\
&The function on the right side of the equation & $F\left(t, \mathbf{x}, u(t, \mathbf{x}),\left(\nabla_{\mathbf{x}} u\right)(t, \mathbf{x}),\left(\operatorname{Hess}_{\mathbf{x}} u\right)(t, \mathbf{x})\right)$\\
\hline
{Stochastic symbol}
&The probability space & $(\Omega, \mathcal{F}, \mathbb{P})$ \\
&The standard Brownian motion & $\mathcal{W}$ \\
&The normal filtration generated via $\mathcal{W}$ & $\mathbb{F}_{t}$(abbr. $\mathbb{F}$) \\
&The $\mathbb{F}$-adapted stochastic process & $\mathcal{X},\mathcal{Y},\mathcal{Z},\Gamma,\mathcal{A}$\\
&The state of the $\mathbb{F}$-adapted stochastic process at time $t$ & $\mathcal{X}_t,\mathcal{Y}_t,\mathcal{Z}_t,\Gamma_t,\mathcal{A}_t$\\
\hline
{Deep learning symbol}
&The approximate function by deep learning & $\mathbf{G}_n^{\theta},\mathbf{A}_n^{\theta},\mathcal{Y}_n^{\theta},\mathcal{Z}_n^{\theta},\tilde{\mathbf{G}}_n^{\theta},\tilde{\mathbf{A}}_n^{\theta}$\\
&The activation function(ReLU function) & $\mathbf{R}_k(\mathbf{x})$\\
&The linear affine function & $\mathbf{M}_{k,l}^{\theta,v}$\\
&The $d_i$ scale neural networks & $\mathbf{G}_{d_i}^{\theta},\mathbf{A}_{d_i}^{\theta}$\\
&The convolution function & $\tilde{\mathbf{M}}_{k,l}^{\theta,v}$\\
&The $i$ channel convolution neural networks & $\tilde{\mathbf{G}}_{i}^{\theta},\tilde{\mathbf{A}}_{i}^{\theta}$\\
&The loss function of training & $\widetilde{\phi}^{m, \mathbf{s}}(\theta, \omega)$\\
&The function of learning rate & $\tilde{\gamma}(m)$(abbr. $\tilde{\gamma}_m$)\\
\hline
{Basic symbol}
&The dimension of $\mathbf{x}$& $d(d\in\mathbb{Z}^+)$ \\
&The range of $t$ & $T(0<T<\infty)$\\
&The number of time discrete points & $N(N\geq 1)$ \\
&The certain point in time & $t_n(0\leq n\leq N-1)$ \\
&The number of parameters in deep learning & $v$ \\
&The number of channels in convolution neural network & $c(c\geq 1)$ \\
\hline
\end{tabular}
\end{center}
\end{table}

\subsection{\normalfont Fully nonlinear second-order PDEs}\label{Se2.1}
\vskip 0.2mm
Let $d \in \mathbb{Z}^+$, $0<T < \infty$, $u=(u(t, \mathbf{x}))_{0\leq t \leq T, \mathbf{x} \in \mathbb{R}^{d}} \in C^{1,3}\left([0, T] \times \mathbb{R}^{d}, \mathbb{R}\right)$, $F \in C([0, T] \times \mathbb{R}^{d} \times \mathbb{R} \times \mathbb{R}^{d} \times \mathbb{R}^{d \times d}, \mathbb{R})$ and $\hat{g} \in C(\mathbb{R}^{d}, \mathbb{R})$ satisfy that $u(T, \mathbf{x})=\hat{g}(\mathbf{x})$ and
\begin{equation}\label{eq2.1}
\begin{split}
    \frac{\partial u}{\partial t}(t, \mathbf{x}) = F\left(t, \mathbf{x}, u(t, \mathbf{x}),\left(\nabla_{\mathbf{x}} u\right)(t, \mathbf{x}),\left(\operatorname{Hess}_{\mathbf{x}} u\right)(t, \mathbf{x})\right),
\end{split}
\end{equation}
for all $t \in[0, T)$ and $\mathbf{x} \in \mathbb{R}^{d}$.

\vskip 0.2mm
Then, the deep-learning 2BSDE approaches can effective approximate the function $u(0,\mathbf{x})\in \mathbb{R}$ with $\mathbf{x}\in \mathbb{R}^{d}$. Note that deep-learning 2BSDE techniques can be easily extended to the case of fully nonlinear second-order parabolic PDEs, but for keeping the symbolic complexity as low as possible, we restrict ourselves to the scalar case in this work (see \eqref{eq2.1}).

\vskip 0.2mm
Furthermore, equation \eqref{eq2.1} is formulated as a terminal value problem. We select the terminal value problem instead of the initial value problem, which is more common in the literature of PDEs. On the one hand, the terminal value problem seems to be more naturally associated with 2BSDEs (see Section \ref{Se2.2}), and on the other hand, the terminal value problem naturally appears in financial engineering applications such as the BSB equation in derivatives pricing (see Section \ref{Se4.2}). Obviously, terminal value problems can be transformed into initial value problems and vice versa, which can be seen in the following Lemma.

\begin{lemma}\normalfont \label{lemma2.1}\cite[Lemma 3.1]{Beck}
  Let $d \in \mathbb{Z}^+$, $0<T < \infty$, $F:[0, T] \times \mathbb{R}^{d} \times \mathbb{R} \times \mathbb{R}^{d} \times \mathbb{R}^{d \times d} \rightarrow \mathbb{R}$ and $\hat{g}: \mathbb{R}^{d} \rightarrow \mathbb{R}$, and assume that $u:[0, T] \times \mathbb{R}^{d} \rightarrow \mathbb{R}$ be a continuous function such that $u(T, \mathbf{x})=\hat{g}(\mathbf{x}), u|_{[0, T) \times \mathbb{R}^{d}} \in C^{1,2}\left([0, T) \times \mathbb{R}^{d}, \mathbb{R}\right)$ and
\begin{equation}\label{eq2.2}
\begin{split}
\frac{\partial u}{\partial t}(t, \mathbf{x})=F\left(t, \mathbf{x}, u(t, \mathbf{x}),\left(\nabla_{\mathbf{x}} u\right)(t, \mathbf{x}),\left(\operatorname{Hess}_{\mathbf{x}} u\right)(t, \mathbf{x})\right),
\end{split}
\end{equation}
for all $(t, \mathbf{x}) \in[0, T) \times \mathbb{R}^{d}$. Assume $\widehat{F}:[0, T] \times \mathbb{R}^{d} \times \mathbb{R} \times \mathbb{R}^{d} \times \mathbb{R}^{d \times d} \rightarrow \mathbb{R}$ and $V:[0, T] \times \mathbb{R}^{d} \rightarrow \mathbb{R}$ be the functions such that $V(t, \mathbf{x})=u(T-t, \mathbf{x})$ and
\begin{equation}\label{eq2.3}
\begin{split}
 \widehat{F}(t, \mathbf{x}, \mathbf{y}, \mathbf{z}, \rho)=-F(T-t, \mathbf{x}, \mathbf{y}, \mathbf{z}, \rho),
\end{split}
\end{equation}
for all $(t, \mathbf{x}, \mathbf{y}, \mathbf{z}, \rho) \in[0, T] \times \mathbb{R}^{d} \times \mathbb{R} \times \mathbb{R}^{d} \times \mathbb{R}^{d \times d}$. Then we get that $V: [0, T] \times \mathbb{R}^{d} \rightarrow \mathbb{R}$ is a continuous function, such that $V(0, x)=\hat{g}(x), V|_{(0, T] \times \mathbb{R}^{d}} \in C^{1,3}\left((0, T] \times \mathbb{R}^{d}, \mathbb{R}\right)$ and
\begin{equation}\label{eq2.4}
\begin{split}
\frac{\partial V}{\partial t}(t, \mathbf{x})=\widehat{F}\left(t, \mathbf{x}, V(t, \mathbf{x}),\left(\nabla_{\mathbf{x}} V\right)(t, \mathbf{x}),\left(\operatorname{Hess}_{\mathbf{x}} V\right)(t, \mathbf{x})\right),
\end{split}
\end{equation}
for all $(t, \mathbf{x}) \in(0, T] \times \mathbb{R}^{d}$.
\end{lemma}
Based on the above discussion, in the following numerical examples, we only consider the terminal problem.
\subsection{\normalfont Combination between fully nonlinear second-order PDEs and 2BSDEs}\label{Se2.2}
\vskip 0.2mm
We apply the deep-learning 2BSDE approaches depend on a combination between fully nonlinear second-order PDEs and 2BSDEs (see the following Lemma \ref{lemma2.2}), from which, It\^{o}' lemma and some suitable assumptions are employed (see \cite{Beck}).

\begin{lemma}\normalfont \label{lemma2.2}
\cite[Lemma 3.1]{Beck}
Assume that $d \in \mathbb{Z}^+$, $0<T < \infty$, and that $u=(u(t, \mathbf{x}))_{t \in[0, T], \mathbf{x} \in \mathbb{R}^{d}} \in C^{1,3}\left([0, T] \times \mathbb{R}^{d}, \mathbb{R}\right), \mu \in C\left(\mathbb{R}^{d}, \mathbb{R}^{d}\right), \sigma \in C\left(\mathbb{R}^{d}, \mathbb{R}^{d \times d}\right), F:[0, T] \times$ $\mathbb{R}^{d} \times \mathbb{R} \times \mathbb{R}^{d} \times \mathbb{R}^{d \times d} \rightarrow \mathbb{R}$, and $\hat{g}: \mathbb{R}^{d} \rightarrow \mathbb{R}$ be functions such that $\nabla_{\mathbf{x}} u \in C^{1,2}\left([0, T] \times \mathbb{R}^{d}, \mathbb{R}^{d}\right), u(T, \mathbf{x})=\hat{g}(\mathbf{x})$ and
\begin{equation}\label{eq2.5}
\begin{split}
\frac{\partial u}{\partial t}(t, \mathbf{x})=F\left(t, \mathbf{x}, u(t, \mathbf{x}),\left(\nabla_{\mathbf{x}} u\right)(t, \mathbf{x}),\left(\operatorname{Hess}_{\mathbf{x}} u\right)(t, \mathbf{x})\right),
\end{split}
\end{equation}
for all $t \in[0, T)$ and $\mathbf{x} \in \mathbb{R}^{d}$. Then, assume that $(\Omega, \mathcal{F}, \mathbb{P})$ is a probability space, that $\mathcal{W}=\left(\mathcal{W}^{(1)}, \ldots, \mathcal{W}^{(d)}\right):[0, T] \times \Omega \rightarrow \mathbb{R}^{d}$ is a standard Brownian motion on $(\Omega, \mathcal{F}, \mathbb{P})$, that $\mathbb{F}=\left(\mathbb{F}_{t}\right)_{t \in[0, T]}$ is the normal filtration on $(\Omega, \mathcal{F}, \mathbb{P})$ generated via $\mathcal{W}$, that $\xi: \Omega \rightarrow \mathbb{R}^{d}$ is a $\mathcal{F}_{0} / \mathcal{B}\left(\mathbb{R}^{d}\right)$-measurable function, and that $\mathcal{X}=$ $\left(\mathcal{X}^{(1)}, \ldots, \mathcal{X}^{(d)}\right):[0, T] \times \Omega \rightarrow \mathbb{R}^{d}$ is an $\mathbb{F}$-adapted stochastic process, with continuous sample paths such that for all $0\leq t \leq T$, it holds $\mathbb{P}$-a.s. that
\begin{equation}\label{eq2.6}
\begin{split}
\mathcal{X}_{t}=\xi+\int_{0}^{t} \mu\left(\mathcal{X}_{s}\right) {\rm d} s+\int_{0}^{t} \sigma\left(\mathcal{X}_{s}\right) {\rm d} \mathcal{W}_{s},
\end{split}
\end{equation}
for all $\varpi \in C^{1,3}\left([0, T] \times \mathbb{R}^{d}, \mathbb{R}\right)$, and let $\mathcal{L} \varpi: [0, T] \times \mathbb{R}^{d} \rightarrow \mathbb{R}$ be the function such that
\begin{equation}\label{eq2.7}
\begin{split}
(\mathcal{L} \varpi)(t, \mathbf{x})=\left(\frac{\partial \varpi}{\partial t}\right)(t, \mathbf{x})+\frac{1}{2} \operatorname{Trace}\left(\sigma(\mathbf{x}) \sigma(\mathbf{x})^{*}\left(\operatorname{Hess}_{\mathbf{x}} \varpi\right)(t, \mathbf{x})\right),
\end{split}
\end{equation}
for all $(t, \mathbf{x}) \in[0, T] \times \mathbb{R}^{d}$, and let $\mathcal{Y}: [0, T] \times \Omega \rightarrow \mathbb{R}$, $\mathcal{Z}=\left(\mathcal{Z}^{(1)}, \ldots, \mathcal{Z}^{(d)}\right): [0, T] \times \Omega \rightarrow \mathbb{R}^{d}, \Gamma=\left(\Gamma^{(i, j)}\right)_{(i, j) \in\{1, \ldots, d\}^{2}}: [0, T] \times$ $\Omega \rightarrow \mathbb{R}^{d \times d}$, and let $\mathcal{A}=\left(\mathcal{A}^{(1)}, \ldots, \mathcal{A}^{(d)}\right):[0, T] \times \Omega \rightarrow \mathbb{R}^{d}$ be the stochastic processes, such that
\begin{equation}\label{eq2.8}
\begin{split}
\mathcal{Y}_{t}=u\left(t, \mathcal{X}_{t}\right), \;\; \mathcal{Z}_{t}=\left(\nabla_{x} u\right)\left(t, X_{t}\right), \;\; \Gamma_{t}=\left(\operatorname{Hess}_{\mathbf{x}} u\right)\left(t, \mathcal{X}_{t}\right), \;\; \mathcal{A}_{t}^{(i)}=\left(\mathcal{L}\left(\frac{\partial u}{\partial \mathbf{x}_{i}}\right)\right)\left(t, \mathcal{X}_{t}\right)
\end{split}
\end{equation}
for all $0\leq t \leq T$ and $i \in\{1,2, \ldots, d\}$. Then, we obtain that $\mathcal{Y}, \mathcal{Z}, \Gamma, \mathcal{A}$ are $\mathbb{F}$-adapted stochastic processes, with continuous sample paths which satisfy that for all $0\leq t \leq T$, it holds $\mathbb{P}$-a.s. that
\begin{equation}\label{eq2.9}
\begin{split}
\mathcal{Y}_{t} =\hat{g}\left(\mathcal{X}_{T}\right) &  -\int_{t}^{T}\left(F\left(s, \mathcal{X}_{s}, \mathcal{Y}_{s}, \mathcal{Z}_{s}, \Gamma_{s}\right)+\frac{1}{2} \operatorname{Trace}\left(\sigma\left(\mathcal{X}_{s}\right) \sigma\left(\mathcal{X}_{s}\right)^{*} \Gamma_{s}\right)\right) {\rm d}s \\
 &-\int_{t}^{T}\left\langle \mathcal{Z}_{s}, {\rm d}\mathcal{X}_{s}\right\rangle_{\mathbb{R}^{d}}
\end{split}
\end{equation}
and
\begin{equation}\label{eq2.10}
\begin{split}
\mathcal{Z}_{t}=\mathcal{Z}_{0}+\int_{0}^{t} \mathcal{A}_{s} {\rm d}s+\int_{0}^{t} \Gamma_{s} {\rm d} \mathcal{X}_{s}.
\end{split}
\end{equation}

\end{lemma}

\subsection{\normalfont Merged construction of PDEs and 2BSDEs}\label{Se2.3}
\vskip 0.2mm
In what follows, we present a merged construction for PDE \eqref{eq2.1} and 2BSDE system.
Let the hypotheses in Lemma \ref{lemma2.2} be satisfied and use the same notations as Lemma \ref{lemma2.2}. Then, one can easily see that for $0\leq \delta_{1}, \delta_{2} \leq T$,
\begin{equation}\label{eq2.11}
    \begin{split}
        \mathcal{X}_{\delta_{2}}=\mathcal{X}_{\delta_{1}}+\int_{\delta_{1}}^{\delta_{2}} \mu\left(\mathcal{X}_{s}\right) {\rm d} s+\int_{\delta_{1}}^{\delta_{2}} \sigma\left(\mathcal{X}_{s}\right) {\rm d} \mathcal{W}_{s},
    \end{split}
\end{equation}
\begin{equation}\label{eq2.12}
\begin{split}
    &\mathcal{Y}_{\delta_{2}}=\mathcal{Y}_{\delta_{1}}+\int_{\delta_{1}}^{\delta_{2}}\left\langle \mathcal{Z}_{s}, d \mathcal{X}_{s}\right\rangle_{\mathbb{R}^{d}} \\
    &+\int_{\delta_{1}}^{\delta_{2}}\left(F\left(s, \mathcal{X}_{s}, \mathcal{Y}_{s}, \mathcal{Z}_{s},\left(\operatorname{Hess}_{\mathbf{x}} u\right)\left(s, \mathcal{X}_{s}\right)\right)+\frac{1}{2} \operatorname{Trace}\left(\sigma\left(\mathcal{X}_{s}\right) \sigma\left(\mathcal{X}_{s}\right)^{*}\left(\operatorname{Hess}_{\mathbf{x}} u\right)\left(s, \mathcal{X}_{s}\right)\right)\right) {\rm d}s
\end{split}
\end{equation}
and
\begin{equation}\label{eq2.13}
\begin{split}
\mathcal{Z}_{\delta_{2}}=\mathcal{Z}_{\delta_{1}}+\int_{\delta_{1}}^{\delta_{2}}\left(\mathcal{L}\left(\nabla_{\mathbf{x}} u\right)\right)\left(s, \mathcal{X}_{s}\right) {\rm d} s+\int_{\delta_{1}}^{\delta_{2}}\left(\operatorname{Hess}_{\mathbf{x}} u\right)\left(s, \mathcal{X}_{s}\right) {\rm d} \mathcal{X}_{s}.
\end{split}
\end{equation}

\section{Approximation of the merged PDEs-2BSDEs system}\label{Se3}

\subsection{\normalfont Forward-discretizations of the merged PDEs-2BSDEs system}\label{Se3.1}
\vskip 0.2mm
Now, we  describe a forward discretization of the merged PDEs-2BSDEs system \eqref{eq2.11}-\eqref{eq2.13}.  Let us consider  positive integer $N\geq 1$ with $t_{0}, t_{1}, \ldots, t_{N} \in [0, T]$, such that
$$
0=t_{0}<t_{1}<t_{2}<\ldots<t_{N}=T,
$$
from which, the max mesh size $\tau:=\max\limits_{0 \leq j \leq N-1}\left(t_{j+1}-t_{j}\right)$ is sufficiently small and we define $\tau_j = t_{j}-t_{j-1}$ for $1 \leq j \leq N$.
\vskip 0.2mm
Notice that, for sufficiently large $N \in \mathbb{Z}^+$, \eqref{eq2.6}-\eqref{eq2.8} and \eqref{eq2.11}-\eqref{eq2.13} indicate that for all $n \in\{0,1, \ldots, N-1\}$, it holds that
\begin{equation}\label{eq3.1}
\begin{split}
\mathcal{X}_{t_{0}}=\mathcal{X}_{0}=\xi, \qquad \mathcal{Y}_{t_{0}}=\mathcal{Y}_{0}=u(0, \xi), \qquad \mathcal{Z}_{t_{0}}=\mathcal{Z}_{0}=\left(\nabla_{\mathbf{x}} u\right)(0, \xi),
\end{split}
\end{equation}
\begin{equation}\label{eq3.2}
\begin{split}
\mathcal{X}_{t_{n+1}} \approx \mathcal{X}_{t_{n}}+\mu\left(\mathcal{X}_{t_{n}}\right)\tau_{n+1}+\sigma\left(\mathcal{X}_{t_{n}}\right)\left(\mathcal{X}_{t_{n+1}}-\mathcal{X}_{t_{n}}\right),
\end{split}
\end{equation}
\begin{equation}\label{eq3.3}
\begin{split}
\mathcal{Y}_{t_{n+1}} \approx \mathcal{Y}_{t_{n}} &+ \Big[ F\big( t_{n}, \mathcal{X}_{t_{n}}, \mathcal{Y}_{t_{n}}, \mathcal{Z}_{t_{n}},(\operatorname{Hess}_{\mathbf{x}} u)(t_{n}, \mathcal{X}_{t_{n}}) \big) \\
&+\frac{1}{2} \operatorname{Trace}(\sigma(\mathcal{X}_{t_{n}}) \sigma(\mathcal{X}_{t_{n}})^{*}(\operatorname{Hess}_{\mathbf{x}} u)(t_{n}, \mathcal{X}_{t_{n}})) \Big]  \tau_{n+1} +\langle \mathcal{Z}_{t_{n}}, \mathcal{X}_{t_{n+1}}-\mathcal{X}_{t_{n}}\rangle_{\mathbb{R}^{d}},
\end{split}
\end{equation}
and
\begin{equation}\label{eq3.4}
\begin{split}
\mathcal{Z}_{t_{n+1}} \approx \mathcal{Z}_{t_{n}}+\left(\mathcal{L}\left(\nabla_{\mathbf{x}} u\right)\right)\left(t_{n}, \mathcal{X}_{t_{n}}\right)\tau_{n+1} +\left(\operatorname{Hess}_{\mathbf{x}} u\right)\left(t_{n}, \mathcal{X}_{t_{n}}\right)\left(\mathcal{X}_{t_{n+1}}-\mathcal{X}_{t_{n}}\right).
\end{split}
\end{equation}
Naturally, we can obtain the semi-discretization approximation of the merged PDEs-2BSDEs system by \eqref{eq3.1}-\eqref{eq3.4}.

\subsection{\normalfont Spatial discretizations based on multiscale deep learning fusion}\label{Se3.2}

\vskip 0.2mm
In the following, for all $0\leq n \leq N-1$ and $\mathbf{x}\in \mathbb{R}^{d}$, we select suitable approximations for functions
$\left(\operatorname{Hess}_{\mathbf{x}} u\right)\left(t_{n}, \mathbf{x}\right) \in \mathbb{R}^{d \times d}$ and
$ \left(\mathcal{L}(\nabla_\mathbf{x} u)\right)\left(t_{n}, \mathbf{x}\right) \in \mathbb{R}^{d} $
given in \eqref{eq3.3}-\eqref{eq3.4} and for the functions $u\left(t_{n}, \mathbf{x}\right) \in \mathbb{R}^{d}$ and $\left(\nabla_{\mathbf{x}} u\right)\left(t_{n}, \mathbf{x}\right) \in \mathbb{R}^{d}$. Precisely, we assume that $\nu \in \mathbb{N} \cap[d+1, \infty)$ for every $\theta \in \mathbb{R}^{\nu}, 0\leq n \leq N$.

\vskip 0.2mm
Assume $\mathbf{G}_{n}^{\theta}: \mathbb{R}^{d} \rightarrow \mathbb{R}^{d \times d}$ and $\mathbf{A}_{n}^{\theta}: \mathbb{R}^{d} \rightarrow \mathbb{R}^{d}$ are continuous functions, and then, for every $\theta=$ $\left(\theta_{1}, \theta_{2}, \ldots, \theta_{\nu}\right) \in \mathbb{R}^{\nu}$, assume $\mathcal{Y}^{\theta}:\{0,1, \ldots, N\} \times \Omega \rightarrow \mathbb{R}$ and $\mathcal{Z}^{\theta}:\{0,1, \ldots, N\} \times \Omega \rightarrow \mathbb{R}^{d}$ be stochastic processes, such that $\mathcal{Y}_{0}^{\theta}=\theta_{1}$, $\mathcal{Z}_{0}^{\theta}=\left(\theta_{2}, \theta_{3}, \ldots, \theta_{d+1}\right)$,
\begin{equation}\label{eq3.5}
\begin{split}
\mathcal{Y}_{n+1}^{\theta}&=\mathcal{Y}_{n}^{\theta}+\left\langle\mathcal{Z}_{n}^{\theta}, \mathcal{X}_{t_{n+1}}-\mathcal{X}_{t_{n}}\right\rangle_{\mathbb{R}^{d}} \\
&+\left(F\left(t_{n}, \mathcal{X}_n, \mathcal{Y}_{n}^{\theta}, \mathcal{Z}_{n}^{\theta}, \mathbf{G}_{n}^{\theta}\left(\mathcal{X}_n\right)\right)+\frac{1}{2} \operatorname{Trace}\left(\mathbf{G}_{n}^{\theta}\left(\mathcal{X}_n\right)\right)\right)\tau_{n+1}
\end{split}
\end{equation}
and
\begin{equation}\label{eq3.6}
\begin{split}
\mathcal{Z}_{n+1}^{\theta}=\mathcal{Z}_{n}^{\theta}+\mathbf{A}_{n}^{\theta}\left(\mathcal{X}_n\right)\tau_{n+1} +\mathbf{G}_{n}^{\theta}\left(\mathcal{X}_n\right)\left(\mathcal{X}_{t_{n+1}}-\mathcal{X}_{t_{n}}\right),
\end{split}
\end{equation}
for $0\leq n \leq N-1$. For all favorable $\theta \in \mathbb{R}^{\nu}$, $\mathbf{x} \in \mathbb{R}^{d}$ and $0\leq n \leq N-1$, we select the suitable approximations that
$\mathcal{Y}_{n}^{\theta} \approx \mathcal{Y}_{t_{n}} $, $ \mathcal{Z}_{n}^{\theta} \approx \mathcal{Z}_{t_{n}}$,
$\mathbf{G}_{n}^{\theta}(\mathbf{x}) \approx\left(\text {Hess}_{\mathbf{x}} u\right)\left(t_{n}, \mathbf{x}\right) $ and
$\mathbf{A}_{n}^{\theta}(\mathbf{x}) \approx \left(\mathcal{L}(\nabla_\mathbf{x} u)\right)\left(t_{n}, \mathbf{x}\right)$.

\vskip 0.2mm
Especially, we regard $\theta_{1}$ and $\left(\theta_{2}, \theta_{3}, \ldots, \theta_{d+1}\right)$ as the suitable approximations of $u(0, \xi)$ and $\left(\nabla_{\mathbf{x}} u\right)(0, \xi)$ with $u(0, \xi) \in \mathbb{R}$ and $\left(\nabla_{\mathbf{x}} u\right)(0, \xi) \in \mathbb{R}^{d}$. Then we can select the functions $\mathbf{G}_{d_i}^{\theta}$ and $\mathbf{A}_{d_i}^{\theta}$ as deep neural networks. In particular, $d_i$ represents the scale of different neural networks, and four scales are selected here. Furthermore, we use the same neural network for different time $n$. That is, the parameters of our network only depend on different scales, independent of time $n$.

\vskip 0.2mm
Assume $\nu \geq\left(2\sum_{i=1}^4 d_i+d+1\right)(d+1)+\sum_{i=1}^4(2d_i+d^2+d)(d_i+1)$. Supposing for all $\theta=\left(\theta_{1}, \ldots, \theta_{\nu}\right) \in \mathbb{R}^{\nu}, \mathbf{x} \in \mathbb{R}^{d}$, we have
$$
\mathbf{G}_{0}^{\theta}(\mathbf{x})=\left(\begin{array}{cccc}
\theta_{d+2} & \theta_{d+3} & \ldots & \theta_{2 d+1} \\
\theta_{2 d+2} & \theta_{2 d+3} & \ldots & \theta_{3 d+1} \\
\vdots & \vdots & \vdots & \vdots \\
\theta_{d^{2}+2} & \theta_{d^{3}+3} & \ldots & \theta_{d^{2}+d+1}
\end{array}\right) \in \mathbb{R}^{d \times d} \text { and } \mathbf{A}_{0}^{\theta}(\mathbf{x})=\left(\begin{array}{c}
\theta_{d^{2}+d+2} \\
\theta_{d^{2}+d+3} \\
\vdots \\
\theta_{d^{2}+2 d+1}
\end{array}\right) \in \mathbb{R}^{d}.
$$
With all $k \in \mathbb{N}$, we let $\mathbf{R}_{k}: \mathbb{R}^{k} \rightarrow \mathbb{R}^{k}$ be the activation function (ReLU), such that
\begin{equation}
\label{eq3.7}
\mathbf{R}_{k}(\mathbf{x})=\Big(\max \left\{\mathbf{x}_{1}, 0\right\}, \ldots, \max \left\{\mathbf{x}_{k}, 0\right\}\Big),
\end{equation}
for every $\mathbf{x}=\left(\mathbf{x}_{1}, \ldots, \mathbf{x}_{k}\right) \in \mathbb{R}^{k}$.  For every $\theta=\left(\theta_{1}, \ldots, \theta_{\nu}\right) \in \mathbb{R}^{\nu}$, $v \in \mathbb{N}_{0}$, $k, l \in \mathbb{N}$ and $v+k(l+1) \leq \nu$, assume $\mathbf{M}_{k, l}^{\theta, v}: \mathbb{R}^{l} \rightarrow \mathbb{R}^{k}$ is the affine linear function, such that
\begin{equation}
\label{eq3.8}
\mathbf{M}_{k, l}^{\theta, v}(\mathbf{x})=\left(\begin{array}{cccc}
\theta_{v+1} & \theta_{v+2} & \ldots & \theta_{v+l} \\
\theta_{v+l+1} & \theta_{v+l+2} & \ldots & \theta_{v+2 l} \\
\theta_{v+2 l+1} & \theta_{v+2 l+2} & \ldots & \theta_{v+3 l} \\
\vdots & \vdots & \vdots & \vdots \\
\theta_{v+(k-1) l+1} & \theta_{v+(k-1) l+2} & \ldots & \theta_{v+k l}
\end{array}\right)\left(\begin{array}{c}
\mathbf{x}_{1} \\
\mathbf{x}_{2} \\
\mathbf{x}_{3} \\
\vdots \\
\mathbf{x}_{l}
\end{array}\right)+\left(\begin{array}{c}
\theta_{v+k l+1} \\
\theta_{v+k l+2} \\
\theta_{v+k l+3} \\
\vdots \\
\theta_{v+k l+k}
\end{array}\right) \\
:= \mathbf{P} \mathbf{x}^T + \mathbf{Q},
\end{equation}
for all $\mathbf{x}=\left(\mathbf{x}_{1}, \ldots, \mathbf{x}_{l}\right)$. For every $\theta \in \mathbb{R}^{\nu}$, $\{d_i|i\in\{1,2,3,4\},d_0=0\},d^{t}=\sum_{i=1}^4 d_i,d^{t^2}=\sum_{i=1}^4 d_i^2+d_i,d^{t^3}=\sum_{i=1}^4 (d_i+d)(d_i+1)$ and $\mathbf{x} \in \mathbb{R}^{d}$, we assume that
\begin{eqnarray}\label{Eq:A1}
\mathbf{A}_{d_i}^{\theta}&=&\mathbf{M}_{d, d_i}^{\theta,(d^{t}+d+1)(d+1)+d^{t^2}+\sum_{1}^{i-1}d(d_i+1)} \circ \mathbf{R}_{d_i} \circ \notag \\
&&\mathbf{M}_{d_i, d_i}^{\theta,(d^{t}+d+1)(d+1)+\sum_{1}^{i-1}d_i^2+d_i} \circ \mathbf{R}_{d_i} \circ \mathbf{M}_{d_i, d}^{\theta,(\sum_{1}^{i-1} d_i+d+1)(d+1)},
\end{eqnarray}
and that
\begin{eqnarray}\label{Eq:G1}
\mathbf{G}_{d_i}^{\theta}&=&\mathbf{M}_{d^{2}, d_i}^{\theta,(2d^t+d+1)(d+1)+d^{t^3}+d^{t^2}+\sum_{1}^{i-1}d^2(d_i+1)} \circ \mathbf{R}_{d_i} \circ \notag \\
&&\mathbf{M}_{d_i, d_i}^{\theta,(2d^t+d+1)(d+1)+d^{t^3}+\sum_{1}^{i-1}d_i^2+d_i} \circ \mathbf{R}_{d_i} \circ \mathbf{M}_{d_i, d}^{\theta,(\sum_{1}^{i-1} d_i+d^t+d+1)(d+1)+d^{t^3}}.
\end{eqnarray}

\begin{figure}[thbp]
\hspace{0.8cm}\includegraphics[width=0.9\textwidth]{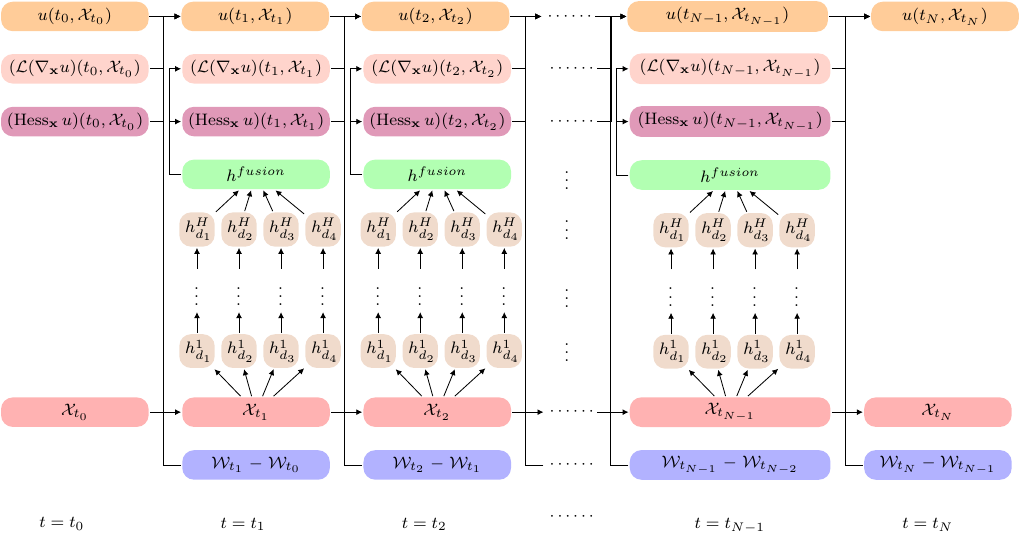}
\caption{\quad Sketch of the architecture of the multiscale deep learning fusion for BSDE.}
\label{TK1}
\end{figure}

\begin{remark}\label{re3.1} \normalfont
In this remark, we illustrate the multiscale deep learning fusion and the specific choice of the $\nu\in\mathbb{N}$ in the above.
\begin{enumerate}
    \item[(i)] Multiscale fusion is mainly reflected in function $\mathbf{A}_{d_i}^{\theta}$ and $\mathbf{G}_{d_i}^{\theta}$. We use deep neural networks of different scales to obtain $\mathbf{A}_{d_i}^{\theta}$ and $\mathbf{G}_{d_i}^{\theta}$, then fuse them to get the final result. In fact, multiscale fusion is to obtain more information in neural network training, thereby improving training results. In addition, if it is assumed that the scales selected each time are the same, then our multiscale fusion is equivalent to a weighted average of multiple experiments. From a probabilistic point of view, the results of multiple experiments are often more accurate and stable than the results of a single experiment.
    \item[(ii)] For the specific choice of the $\nu$, the choice of $\nu$ is mainly divided into three parts. On the one hand, it is employed to approximate the variables we need, which includes the real number $u(0,\xi)\in\mathbb{R}$, the $(1\times d)$ matrix $(\nabla_{\mathbf{x}} u)(0, \xi)$, the $(d\times d)$ matrix $\mathbf{G}_0^{\theta}$ and $(d\times 1)$ vector $\mathbf{A}_0^{\theta}$. So, we have $\nu\geq(d+1)(d+1)$. On the other hand, the remaining two parts are related to neural networks, the first part is about $\mathbf{A}_{d_i}^{\theta}$, and the last part is about $\mathbf{G}_{d_i}^{\theta}$.
    \item[(iii)] For the $\mathbf{A}_{d_i}^{\theta}$, in each of the employed $d_i$ neural network we use $d_i(d+1)$ components of $\theta$ to describe the affine linear function from the $d$-dimensional first layer (input layer) to the $d_i$-dimensional second layer (includes a $d_i\times d$ matrix and a $d_i$ vector, see \eqref{eq3.8}). Next, we use $d_i(d_i+1)$ to describe the $d_i$-dimensional second layer to the $d_i$-dimensional third layer. Finally, the $d(d_i+1)$ is used in the $d_i$-dimensional third layer to the $d_i$-dimensional fourth layer (output layer). For the $\mathbf{G}_{d_i}^{\theta}$, the few layers are basically the same as $\mathbf{A}_{d_i}^{\theta}$, the only difference is that the $d^2(d_i+1)$ is used in the $d_i$-dimensional third layer to the $d_i$-dimensional fourth layer. Therefore, combining the above analysis, we have
    \begin{eqnarray*}
    \label{remark1}
    v&\geq&(d+1)(d+1)+2\sum_{i=1}^4 d_i(d+1)+2\sum_{i=1}^4 d_i(d_i+1)+\sum_{i=1}^4 d(d_i+1)+\sum_{i=1}^4d^2(d_i+1)\\
    &=&(2\sum_{i=1}^4 d_i+d+1)(d+1)+\sum_{i=1}^4(2d_i+d^2+d)(d_i+1).
    \end{eqnarray*}
    \item[(iv)] We also depict the sketch of the architecture of multiscale deep learning fusion, see Figure \ref{TK1}. In Figure \ref{TK1}, when $t=t_0$, we first give the initial values $\mathcal{X}_{t_0},(\mathcal{L}(\nabla_x u)(t_0,\mathcal{X}_{t_0}),(\mathrm{Hess}_x u)(t_0,\mathcal{X}_{t_0})$. Then use the initial value to calculate the variables of $t=t_i$ in turn ($1\geq i\geq N-1$), which $h_{d_i}^H$ represents the $H$ layer of the neural network at the $d_i$ scale. As can be seen from the figure, each $\mathcal{X}_{t_i}$ is trained by neural networks of four scales, and finally fused to obtain $h^{fusion}$. Note that for each time $t=t_i$, we use the same neural network, which is continuously updated as time changes. In addition, $(\mathcal{L}(\nabla_{\mathbf{x}} u))(t_i,\mathcal{X}_{t_i})$ and $(\mathrm{Hess}_{\mathbf{x}}u)(t_i,\mathcal{X}_{t_i})$ are approximated separately using two networks. In Figure \ref{TK1}, it is not subdivided for the sake of brevity.

\end{enumerate}
\end{remark}

\subsection{\normalfont Spatial discretizations based on convolutional neural network}\label{Se3.3}

\vskip 0.2mm
Here, with $0\leq n \leq N-1$ and $\mathbf{x}\in \mathbb{R}^{d}$, based on convolutional neural network, we still choose the suitable approximations for functions
$\left(\operatorname{Hess}_{\mathbf{x}} u\right)\left(t_{n}, \mathbf{x}\right) \in \mathbb{R}^{d \times d}$,
$ \left(\mathcal{L}(\nabla_\mathbf{x} u)\right)\left(t_{n}, \mathbf{x}\right) \in \mathbb{R}^{d} $,  $u\left(t_{n}, \mathbf{x}\right) \in \mathbb{R}^{d}$ and $\left(\nabla_{\mathbf{x}} u\right)\left(t_{n}, \mathbf{x}\right) \in \mathbb{R}^{d}$. Then let $\nu \in \mathbb{N} \cap[d+1, \infty)$ and $\theta$ is assumed as Subsection \ref{Se3.2} with $0\leq n \leq N$.

\vskip 0.2mm
Suppose that $\widetilde{\mathbf{G}}_{n}^{\theta}: \mathbb{R}^{d} \rightarrow \mathbb{R}^{d \times d}$ and $\widetilde{\mathbf{A}}_{n}^{\theta}: \mathbb{R}^{d} \rightarrow \mathbb{R}^{d}$ are continuous functions. For every $\theta=$ $\left(\theta_{1}, \theta_{2}, \ldots, \theta_{\nu}\right) \in \mathbb{R}^{\nu}$, assume $\mathcal{Y}^{\theta}$ and $\mathcal{Z}^{\theta}$ be denoted as before, which satisfy $\mathcal{Y}_{0}^{\theta}=\theta_{1}$, $\mathcal{Z}_{0}^{\theta}=\left(\theta_{2}, \theta_{3}, \ldots, \theta_{d+1}\right)$,
\begin{equation}\label{eq3.9}
\begin{split}
\mathcal{Y}_{n+1}^{\theta}&=\mathcal{Y}_{n}^{\theta}+\left\langle\mathcal{Z}_{n}^{\theta}, \mathcal{X}_{t_{n+1}}-\mathcal{X}_{t_{n}}\right\rangle_{\mathbb{R}^{d}} \\
&+\left(F\left(t_{n}, \mathcal{X}_n, \mathcal{Y}_{n}^{\theta}, \mathcal{Z}_{n}^{\theta}, \widetilde{\mathbf{G}}_{n}^{\theta}\left(\mathcal{X}_n\right)\right)+\frac{1}{2} \operatorname{Trace}\left(\widetilde{\mathbf{G}}_{n}^{\theta}\left(\mathcal{X}_n\right)\right)\right)\tau_{n+1}
\end{split}
\end{equation}
and that
\begin{equation}\label{eq3.10}
\begin{split}
\mathcal{Z}_{n+1}^{\theta}=\mathcal{Z}_{n}^{\theta}+\widetilde{\mathbf{A}}_{n}^{\theta}\left(\mathcal{X}_n\right)\tau_{n+1} +\widetilde{\mathbf{G}}_{n}^{\theta}\left(\mathcal{X}_n\right)\left(\mathcal{X}_{t_{n+1}}-\mathcal{X}_{t_{n}}\right),
\end{split}
\end{equation}
for $0\leq n \leq N-1$. Then, we can choose suitable approximations that
$\mathcal{Y}_{n}^{\theta} \approx \mathcal{Y}_{t_{n}} $, $ \mathcal{Z}_{n}^{\theta} \approx \mathcal{Z}_{t_{n}}$,
$\widetilde{\mathbf{G}}_{n}^{\theta}(\mathbf{x}) \approx\left(\text {Hess}_{\mathbf{x}} u\right)\left(t_{n}, \mathbf{x}\right) $ and
$\widetilde{\mathbf{A}}_{n}^{\theta}(\mathbf{x}) \approx \left(\mathcal{L}(\nabla_\mathbf{x} u)\right)\left(t_{n}, \mathbf{x}\right)$, in view of convolutional neural network. In addition, we consider $\theta_{1}$ and $\left(\theta_{2}, \theta_{3}, \ldots, \theta_{d+1}\right)$ as the affable approximations of $u(0, \xi)$ and $\left(\nabla_{\mathbf{x}} u\right)(0, \xi)$. Also, we can choose functions $\widetilde{\mathbf{G}}_{n}^{\theta}$ and $\widetilde{\mathbf{A}}_{n}^{\theta}$ as deep convolutional neural networks with $0 \leq n \leq N-1$.

\vskip 0.2mm
Similarly, as in Subsection \ref{Se3.2}, we use the same neural network for $\widetilde{\mathbf{A}}_{n}^{\theta},\widetilde{\mathbf{G}}_{n}^{\theta}, \forall n$. The difference is that we introduce the channel $c$ of the convolution kernel. Therefore, we use the new notation $\widetilde{\mathbf{A}}_{i}^{\theta},\widetilde{\mathbf{G}}_{i}^{\theta}, 1\leq i\leq c,i\in\mathbb{N}$. Suppose $\nu \geq[(4 c+4)d+d^2+1](d+1)$ and for every $\theta=\left(\theta_{1}, \ldots, \theta_{\nu}\right) \in \mathbb{R}^{\nu}, \mathbf{x} \in \mathbb{R}^{d}$, we yield that $\widetilde{\mathbf{G}}_{0}^{\theta}(\mathbf{x})=\mathbf{G}_{0}^{\theta}(\mathbf{x})$ and $\widetilde{\mathbf{A}}_{0}^{\theta}(\mathbf{x})=\mathbf{A}_{0}^{\theta}(\mathbf{x})$. Assume $k \in \mathbb{N}$, and we let the activation function (ReLU) $\mathbf{R}_{k}(\mathbf{x})$ be given in \eqref{eq3.7} for every $\mathbf{x}=\left(\mathbf{x}_{1}, \ldots, \mathbf{x}_{k}\right) \in \mathbb{R}^{k}$.  For every $\theta=\left(\theta_{1}, \ldots, \theta_{\nu}\right) \in \mathbb{R}^{\nu}$, $v \in \mathbb{N}_{0}$, $k, l \in \mathbb{N}$ and $v+k(l+1) \leq \nu$, suppose that $\widetilde{\mathbf{M}}_{k, l}^{\theta, v}: \mathbb{R}^{l} \rightarrow \mathbb{R}^{k}$ satisfies that
\begin{equation}\label{eq3.11}
\begin{split}
\widetilde{\mathbf{M}}_{k, l}^{\theta, v}(\mathbf{Z})=\mathbf{P}\otimes \mathbf{Z}+\mathbf{Q},
\end{split}
\end{equation}
where the notation $\otimes$ represents the convolution rule, the matrix
$$
\mathbf{Z}=\left(\begin{array}{cccc}
\mathbf{x}_{1}  & \mathbf{x}_{\sqrt{k}+1} & \ldots & \mathbf{x}_{k-\sqrt{k}+1} \\
\mathbf{x}_{2}  & \mathbf{x}_{\sqrt{k}+2} & \ldots & \mathbf{x}_{k-\sqrt{k}+2} \\
\mathbf{x}_{3}  & \mathbf{x}_{\sqrt{k}+3} & \ldots & \mathbf{x}_{k-\sqrt{k}+3} \\
\vdots & \vdots & \vdots & \vdots \\
\mathbf{x}_{\sqrt{k}} & \mathbf{x}_{2\sqrt{k}} & \ldots & \mathbf{x}_{k}
\end{array}\right),
$$
and $\mathbf{P}$, $\mathbf{Q}$ are presented in \eqref{eq3.8}.

\vskip 0.2mm
For all $\theta \in \mathbb{R}^{\nu}$, $1 \leq i \leq c$ and $\mathbf{x} \in \mathbb{R}^{d}$, suppose that
\begin{equation}\label{Eq:A2}
\widetilde{\mathbf{A}}_{i}^{\theta}=\mathbf{M}_{d,d}^{\theta,[(2 c+2) d+1](d+1)}\circ\mathbf{Re}\left(\mathbf{R}_{d} \circ \widetilde{\mathbf{M}}_{d, d}^{\theta,[(2 c+1) d+1](d+1)} \circ \mathbf{R}_{d} \circ \widetilde{\mathbf{M}}_{d, d}^{\theta,[(c+i) d+1](d+1)} \circ \mathbf{R}_{d} \circ \widetilde{\mathbf{M}}_{d, d}^{\theta,(i d+1)(d+1)}\right),
\end{equation}
and that
\begin{equation}\label{Eq:G2}
\widetilde{\mathbf{G}}_{i}^{\theta}=\mathbf{M}_{d^2,d}^{\theta,[(4 c+4) d+1](d+1)}\circ\mathbf{Re}\left(\widetilde{\mathbf{M}}_{d^{2}, d}^{\theta,[(4 c+3)d+1](d+1)} \circ \mathbf{R}_{d} \circ \widetilde{\mathbf{M}}_{d, d}^{\theta,[(3c+2+i) d+1](d+1)} \circ \mathbf{R}_{d} \circ \widetilde{\mathbf{M}}_{d, d}^{\theta,[(2c+2+i) d+1](d+1)}\right),
\end{equation}
in which $\mathbf{Re}(\cdot)$ denotes the operation to pull the matrix $\mathbf{Z}$ back into the vector $\mathbf{x}$.

\begin{figure}[thbp]
\hspace{0.6cm}\includegraphics[width=0.9\textwidth]{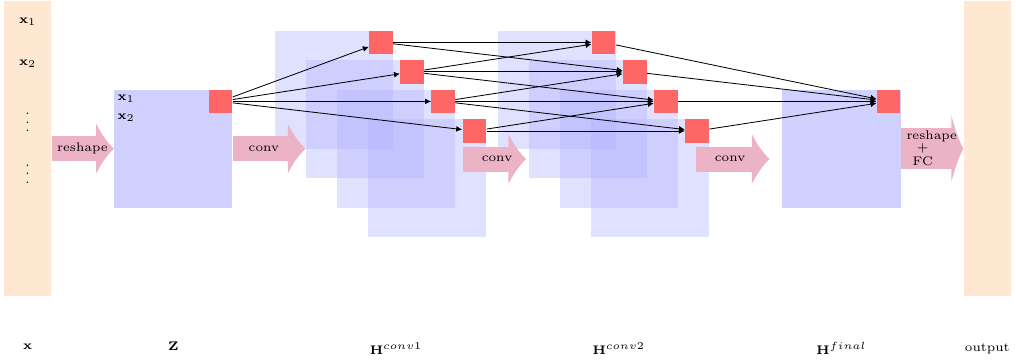}
\caption{\quad The rough schematic diagram of convolutional neural network.}
\label{TK2}
\end{figure}

\begin{remark}\label{re3.2}\normalfont
In this remark, we describe some details in convolutional neural networks.
\begin{enumerate}
    \item[(i)] We used three convolutional layers and one linear layer. In the convolution layer, we use a convolution kernel of $3\times3$, and the stride and padding are both set to $1$ by default. Therefore, the matrix size does not change after each convolution. In the first two convolutional layers, we set the number of channels to $32$, and in the last convolutional layer, set the number of channels to $1$. For the linear layer, we first pull the output of the convolutional layer into vector, then employ the linear transformation in Subsection \ref{Se3.1}.
    \item[(ii)] For the specific choice of the $\nu$, the basic calculation idea is consistent with Subsection \ref{Se3.1}. In first stage, we have $\nu\geq(d+1)(d+1)$ as same as Subsection \ref{Se3.1}. In second stage, the first two convolutional layers are $2c\cdot d(d+1)$, the final convolutional layer is $1\cdot d(d+1)$ and the linear layer is $d(d+1)$ for $\tilde{\mathbf{A}}^{\theta}$. In third stage, for $\tilde{\mathbf{G}}^{\theta}$, except that the linear layer is $d^2(d+1)$, the others are the same as $\tilde{\mathbf{A}}^{\theta}$. We give a specific calculation formula here. For more specific information, please refer to Subsection \ref{Se3.1}.
    \begin{eqnarray*}
    \label{remark2}
    v&\geq&(d+1)(d+1)+2c\cdot d(d+1)+2c\cdot d(d+1)+2\cdot d(d+1)+d(d+1)+d^2(d+1)\\
    &=&[(4 c+4)d+d^2+1](d+1).
    \end{eqnarray*}
    \item[(iii)] Figure \ref{TK2} depicts the rough schematic diagram of convolutional neural network. In fact, other processing processes are similar to Figure \ref{TK1}. For simplicity, we only draw the process of the convolutional neural network here.  As seen in  Figure \ref{TK2}, $\mathbf{x}$ has to undergo a ``reshape'' operation to become $\mathbf{Z}$ before it can be input into the network. As can be seen from the figure, $\mathbf{Z}$ is subjected to a ``conv'' operation to obtain matrix $\mathbf{H}^{conv1}$ of multiple channels. For brevity, only 4 channels are drawn on the graph, there should actually be 32 channels. Note that, in the last layer of convolution $\mathbf{H}^{final}$, we turn the multiple channels back into a single channel. In addition, the ``reshape+FC'' operation means that the matrix is first converted into vector by the ``reshape'' operation. Then ``FC'' is used to perform the operation. Here ``FC'' is the linear transformation in Subsection \ref{Se3.2}.

\end{enumerate}
\end{remark}

\subsection{\normalfont Optimization algorithms} \label{Se3.4}

\vskip 0.2mm
Here, we give the proposed optimization algorithms. First, we present the following lemma (see \cite[Framework 3.2]{Beck}).

\begin{lemma}\normalfont\cite{Beck}\label{lemma3.1}
  Let $T , N, d, \varrho, \varsigma, \nu $ be defined as before. Let $F: [0, T] \times \mathbb{R}^{d} \times$ $\mathbb{R} \times \mathbb{R}^{d} \times \mathbb{R}^{d \times d} \rightarrow \mathbb{R}$ and $\hat{g}: \mathbb{R}^{d} \rightarrow \mathbb{R}$ be functions, and $\left(\Omega, \mathcal{F}, \mathbb{P},\left(\mathbb{F}_{t}\right)_{t \in[0, T]}\right)$ be defined as before. Assume for every $\theta \in \mathbb{R}^{\nu}$ let $\mathbb{U}^{\theta}: \mathbb{R}^{d} \rightarrow \mathbb{R}$ and $\mathbb{Z}^{\theta}: \mathbb{R}^{d} \rightarrow \mathbb{R}^{d}$ be functions and for every $m \in \mathbb{N}_{0}$, $j \in \mathbb{N}$ let $\mathcal{X}^{m, j}:\{0,1, \ldots, N\} \times \Omega \rightarrow \mathbb{R}^{d}$ be a stochastic process such that $\mathcal{X}_{0}^{m, j}=\xi^{m, j}$ and
\begin{equation*}
\begin{split}
\mathcal{X}_{n+1}^{m, j}=\mathcal{H}\left( t_{n}, t_{n+1}, \mathcal{X}_{n}^{m, j}, \mathcal{W}_{t_{n+1}}^{m, j} - \mathcal{W}_{t_{n}}^{m, j} \right),
\end{split}
\end{equation*}
for all $0\leq n \leq N-1$. Then, for every $\theta \in \mathbb{R}^{\nu}, j \in \mathbb{N}, \mathrm{s} \in \mathbb{R}^{\varsigma}, n \in\{0,1, \ldots, N-1\}$, assume $\mathbb{G}_{n}^{\theta, j, \mathrm{~s}}:\left(\mathbb{R}^{d}\right)^{\mathbb{N}_{0}} \rightarrow \mathbb{R}^{d \times d}$ and $\mathbb{A}_{n}^{\theta, j, \mathbf{s}}:\left(\mathbb{R}^{d}\right)^{\mathbb{N}_{0}} \rightarrow \mathbb{R}^{d}$ are functions. Besides, for every $\theta \in \mathbb{R}^{\nu}, m \in \mathbb{N}_{0}, j \in \mathbb{N}, \mathbf{s} \in \mathbb{R}^{\varsigma}$, we suppose that $\mathcal{Y}^{\theta, m, j, \mathrm{~s}}:\{0,1, \ldots, N\} \times \Omega \rightarrow \mathbb{R}$ and $\mathcal{Z}^{\theta, m, j, \mathrm{~s}}:\{0,1, \ldots, N\} \times \Omega \rightarrow \mathbb{R}^{d}$ be stochastic processes such that
\begin{equation*}
\begin{split}
\mathcal{Y}_{0}^{\theta, m, j, \mathbf{s}}=\mathbb{U}^{\theta}\left(\xi^{m, j}\right),  \quad  \mathcal{Z}_{0}^{\theta, m, j, \mathbf{s}}=\mathbb{Z}^{\theta}\left(\xi^{m, j}\right),
\end{split}
\end{equation*}
and
\begin{equation*}
\begin{split}
\mathcal{Y}_{n+1}^{\theta, m, j,\mathbf{s}}& =\mathcal{Y}_{n}^{\theta, m, j,\mathbf{s}}+\tau_{n+1}\Bigl[ \frac{1}{2} \operatorname{Trace}(\sigma(\mathcal{X}_{n}^{m, j}) \sigma(\mathcal{X}_{n}^{m, j})^{*} \mathbf{G}_{n}^{\theta, j, \mathbf{s}}((\mathcal{X}_{n}^{m, i})_{i \in \mathbb{N}}))
\\
& + F(t_{n}, \mathcal{X}_{n}^{m, j}, \mathcal{Y}_{n}^{\theta, m, j}, \mathcal{Z}_{n}^{\theta, m, j, \mathbf{s}}, \mathbf{G}_{n}^{\theta, j, \mathbf{s}}((\mathcal{X}_{n}^{m, i})_{i \in \mathbb{N}})) \Bigl] + \langle\mathcal{Z}_{n}^{\theta, m, j, \mathbf{s}}, \mathcal{X}_{n+1}^{m, j}-\mathcal{X}_{n}^{m, j}\rangle_{\mathbb{R}^{d}}
\end{split}
\end{equation*}
and that
\begin{equation*}
\begin{split}
\mathcal{Z}_{n+1}^{\theta, m, j, \mathbf{s}}=\mathcal{Z}_{n}^{\theta, m, j, \mathbf{s}}+\mathbf{A}_{n}^{\theta, j, \mathbf{s}}\left(\left(\mathcal{X}_{n}^{m, i}\right)_{i \in \mathbb{N}}\right)\tau_{n+1}+\mathbf{G}_{n}^{\theta, j, \mathbf{s}}\left(\left(\mathcal{X}_{n}^{m, i}\right)_{i \in \mathbb{N}}\right)\left(\mathcal{X}_{n+1}^{m, j}-\mathcal{X}_{n}^{m, j}\right).
\end{split}
\end{equation*}
Assume $\left(\mathbf{J}_{m}\right)_{m \in \mathbb{N}_{0}} \subseteq \mathbb{N}$ is a sequence. For every $m \in \mathbb{N}_{0}, \mathbf{s} \in \mathbb{R}^{\varsigma}$, we let $\widetilde{\phi}^{m, \mathbf{s}}: \mathbb{R}^{\nu} \times \Omega \rightarrow \mathbb{R}$ be the function, such that
\begin{equation}\label{Eq:loss}
\begin{split}
\widetilde{\phi}^{m, \mathbf{s}}(\theta, \omega)=\frac{1}{\mathbf{J}_{m}} \sum_{j=1}^{\mathbf{J}_{m}}\left|\mathcal{Y}_{N}^{\theta, m, j, \mathbf{s}}(\omega)-\hat{g}\left(\mathcal{X}_{N}^{m, j}(\omega)\right)\right|^{2}
\end{split}
\end{equation}
for all $(\theta, \omega) \in \mathbb{R}^{\nu} \times \Omega$. Then for every $m \in \mathbb{N}_{0}, \mathrm{~s} \in \mathbb{R}^{\varsigma}$, suppose $\widetilde{\Phi}^{m, \mathrm{~s}}: \mathbb{R}^{\nu} \times \Omega \rightarrow \mathbb{R}^{\nu}$ is a function which satisfies for all $\omega \in \Omega$, $\theta \in\left\{\zeta \in \mathbb{R}^{\nu}: \widetilde{\phi}^{m, \mathbf{s}}(\cdot, \omega): \mathbb{R}^{\nu} \rightarrow \mathbb{R}\right.$ is differentiable at $\left.\zeta\right\}$ that
\begin{equation*}
\begin{split}
\widetilde{\Phi}^{m, \mathbf{s}}(\theta, \omega)=\left(\nabla_{\theta} \widetilde{\phi}^{m, \mathbf{s}}\right)(\theta, \omega),
\end{split}
\end{equation*}
and suppose that $\mathcal{S}: \mathbb{R}^{\varsigma} \times \mathbb{R}^{\nu} \times\left(\mathbb{R}^{d}\right)^{\{0,1, \ldots, N-1\} \times \mathbb{N}} \rightarrow \mathbb{R}^{\varsigma}$ is a function, and for every $m \in \mathbb{N}_{0}$, we let $\widetilde{\psi}_{m}: \mathbb{R}^{\varrho} \rightarrow \mathbb{R}^{\nu}$ and $\widetilde{\Psi}_{m}: \mathbb{R}^{\varrho} \times \mathbb{R}^{\nu} \rightarrow \mathbb{R}^{\varrho}$ be functions. For all $m \in \mathbb{N}_{0}$, we let $\Theta: \mathbb{N}_{0} \times \Omega \rightarrow \mathbb{R}^{\nu}$, $\mathbb{S}: \mathbb{N}_{0} \times \Omega \rightarrow \mathbb{R}^{\varsigma}$, and $\widetilde{\Xi}: \mathbb{N}_{0} \times \Omega \rightarrow \mathbb{R}^{\varrho}$ be stochastic processes, which satisfy that
\begin{equation}\label{eq3.60}
\begin{split}
\mathbb{S}_{m+1}=\mathcal{S}\left(\mathbb{S}_{m}, \Theta_{m},\left(\mathcal{X}_{n}^{m, i}\right)_{(n, i) \in\{0,1, \ldots, N-1\} \times \mathbb{N}}\right),
\end{split}
\end{equation}
and that
\begin{equation}\label{eq3.61}
\begin{split}
\widetilde{\Xi}_{m+1}=\widetilde{\Psi}_{m}\left(\widetilde{\Xi}_{m}, \widetilde{\Phi}^{m, \mathbb{S}_{m+1}}\left(\Theta_{m}\right)\right), \qquad \Theta_{m+1}=\Theta_{m}-\widetilde{\psi}_{m}\left(\widetilde{\Xi}_{m+1}\right).
\end{split}
\end{equation}
\end{lemma}

\vskip 0.05in
Below, we present several special choices for functions $\widetilde{\psi}_{m}, \widetilde{\Psi}_{m}, m \in \mathbb{N}$, given in \eqref{eq3.61}. Based on that, we present the following optimization algorithms.

\vskip 0.2mm
$(i)$ \textbf{Stochastic gradient descent (SGD) method}. Provided the setting in Lemma \ref{lemma3.1}, let notations $\left(\widetilde{\gamma}_{m}\right)_{m \in \mathbb{N}} \subseteq(0, \infty)$, and suppose for all $m \in \mathbb{N}$, $\mathbf{x} \in \mathbb{R}^{\varrho},\left(\varphi_{j}\right)_{j \in \mathbb{N}} \in\left(\mathbb{R}^{\rho}\right)^{\mathbb{N}}$ that
\begin{equation*}
\begin{split}
\varrho=\rho, \qquad \widetilde{\Psi}_{m}\left(\mathbf{x},\left(\varphi_{j}\right)_{j \in \mathbb{N}}\right)=\varphi_{1}, \qquad \widetilde{\psi}_{m}(x)=\widetilde{\gamma}_{m} \mathbf{x},
\end{split}
\end{equation*}
and then it holds that
\begin{equation*}
\begin{split}
\Theta_{m}=\Theta_{m-1}-\widetilde{\gamma}_{m} \widetilde{\Phi}^{m-1}\left(\Theta_{m-1}\right),
\end{split}
\end{equation*}
for all $m \in \mathbb{N}$.

\vskip 0.2mm
$(ii)$ \textbf{Adaptive Moment Estimation (Adam) with mini-batches \cite{57}.} Here, we use Adam optimizer with the deep-learning 2BSDE solver. Provided the setting in Lemma \ref{lemma3.1}, suppose that $\varrho=2 \rho$, and assume $\operatorname{Pow}_{\hat{r}}: \mathbb{R}^{\rho} \rightarrow$ $\mathbb{R}^{\rho}, 0<\hat{r} <\infty$ is the functions satisfying that
\begin{equation*}
\begin{split}
\operatorname{Pow}_{\hat{r}}(x)=\left(\left|\mathbf{x}_{1}\right|^{\hat{r}}, \ldots,\left|\mathbf{x}_{\rho}\right|^{\hat{r}}\right),
\end{split}
\end{equation*}
for all $0<\hat{r} <\infty$ and $\mathbf{x}=\left(\mathbf{x}_{1}, \ldots, \mathbf{x}_{\rho}\right) \in \mathbb{R}^{\rho}$.

\vskip 0.2mm
Let $0<\varepsilon < \infty, \left(\widetilde{\gamma}_{m}\right)_{m \in \mathbb{N}} \subseteq(0, \infty),\left(\mathbf{J}_{m}\right)_{m \in \mathbb{N}_{0}} \subseteq \mathbb{N}$ and $0<\widehat{\mathbb{X}}, \widehat{\mathbb{Y}} <1$, and assume that $\widehat{\mathbf{m}}, \widehat{\mathbb{M}}: \mathbb{N}_{0} \times \Omega \rightarrow \mathbb{R}^{\rho}$ are the stochastic processes which satisfy for all $m \in \mathbb{N}_{0}$ that $\widetilde{\Xi}_{m}=\left(\widehat{\mathbf{m}}_{m}, \widehat{\mathbb{M}}_{m}\right)$, and suppose that
\begin{equation*}
\begin{split}
\widetilde{\Psi}_{m}\left(\mathbf{x}, \mathbf{y},\left(\varphi_{j}\right)_{j \in \mathbb{N}}\right)=\left[\widehat{\mathbb{X}} \mathbf{x}+(1-\widehat{\mathbb{X}})\left(\frac{1}{\mathbf{J}_{m}} \sum_{j=1}^{\mathbf{J}_{m}} \varphi_{j}\right), \widehat{\mathbb{Y}} \mathbf{y}+(1-\widehat{\mathbb{Y}}) \operatorname{Pow}_{2}\left(\frac{1}{\mathbf{J}_{m}} \sum_{j=1}^{\mathbf{J}_{m}} \varphi_{j}\right)\right]
\end{split}
\end{equation*}
and
\begin{equation*}
\begin{split}
\widetilde{\psi}_{m}(\mathbf{x}, \mathbf{y})=\left[\varepsilon+\operatorname{Pow}_{\frac{1}{2}}\left(\mathbf{y}\right)\right]^{-1} \widetilde{\gamma}_{m} \mathbf{x},
\end{split}
\end{equation*}
for all $m \in \mathbb{N}, \mathbf{x}, \mathbf{y} \in \mathbb{R}^{\rho},\left(\varphi_{j}\right)_{j \in \mathbb{N}} \in\left(\mathbb{R}^{\rho}\right)^{\mathbb{N}}$. Then for all $m \in \mathbb{N}$, we have
\begin{equation*}
\begin{split}
\widehat{\mathbb{M}}_{m} =\widehat{\mathbb{Y}} \widehat{\mathbb{M}}_{m-1}+(1-\widehat{\mathbb{Y}}) \operatorname{Pow}_{2}\left(\frac{1}{\mathbf{J}_{m}} \sum_{j=1}^{\mathbf{J}_{m}} \widetilde{\Phi}_{\mathbb{S}_{m}}^{m-1, j}\left(\Theta_{m-1}\right)\right),
\end{split}
\end{equation*}
\begin{equation*}
\begin{split}
\widehat{\mathbf{m}}_{m} =\widehat{\mathbb{X}} \widehat{\mathbf{m}}_{m-1}+(1-\widehat{\mathbb{X}})\left(\frac{1}{\mathbf{J}_{m}}\sum_{j=1}^{\mathbf{J}_{m}} \widetilde{\Phi}_{\mathbb{S}_{m}}^{m-1, j}\left(\Theta_{m-1}\right)\right),
\end{split}
\end{equation*}
and the final update formula is
\begin{equation*}
\begin{split}
\Theta_{m} =\Theta_{m-1}-\left[\varepsilon+\operatorname{Pow}_{\frac{1}{2}}\left(\widehat{\mathbb{M}}_{m}\right)\right]^{-1} \widetilde{\gamma}_{m} \widehat{\mathbf{m}}_{m}.
\end{split}
\end{equation*}

Finally, we summarize the proposed approximation method in Algorithm \ref{AL1}.

\begin{algorithm}
\setcounter{algorithm}{0}
\caption{Approximation algorithm using multi-scale deep learning fusion (or convolutional neural network). }
\label{AL1}
\begin{algorithmic}[1]
\REQUIRE the functions $F(t,\mbx,u(t,\mbx),(\nabla_{\mbx} u)(t,\mbx),(\operatorname{Hess}_{\mbx} u)(t,\mbx))$ and $\hat{g}(\mbx)$.
\ENSURE $u(0,\mbx)$.
\STATE Initialize $\mathcal{X}_{t_0},\mathcal{Y}_{t_0},\mathcal{Z}_{t_0},\mathbf{G},\mathbf{A}$.
\FOR{$t=t_0:t_N$ (each discrete time point)}
\STATE Updated the $\mathcal{X}_t,\mathcal{Y}_t$ by using \eqref{eq3.2} and \eqref{eq3.3}.
\STATE Compute the $\mathbf{A}$ and $\mathbf{G}$ by using \eqref{Eq:A1} and \eqref{Eq:G1} (or \eqref{Eq:A2} and \eqref{Eq:G2}).
\ENDFOR
\WHILE{not up to total training steps}
\STATE Compute the loss function $\widetilde{\phi}^{m, \mathbf{s}}(\theta, \omega)$ by using \eqref{Eq:loss}.
\STATE Apply SGD or Adam algorithm to optimization the loss function $\widetilde{\phi}^{m, \mathbf{s}}(\theta, \omega)$.
\STATE Updated the neural networks (or convolutional neural network) paramters by back propagation.
\ENDWHILE
\IF{the training is completed}
\STATE Obtain the function value $\mathcal{Y}_{t_0}$.
\RETURN $u(0,\mbx)=\mathcal{Y}_{t_0}$.
\ENDIF
\end{algorithmic}
\end{algorithm}

\section{Numerical results and discussion}\label{Se4}
This section employs the multiscale deep learning fusion and CNNs to approximately solve several stochastic PDEs, which mainly include the AC, HJB and BSB equations. Specifically, in Subsection \ref{Se4.1}, we first employ multiscale deep learning to solve the 20-dimensional AC equation and compare with the method of Beck et al. \cite{Beck}, and use the CNNs to obtain numerical solutions of the higher-dimensional AC equation. Then the numerical experiments in 256 and 400 dimensions are given, respectively. Subsections \ref{Se4.2} and \ref{Se4.3} also deal with the HJB and BSB equations, respectively, and the only difference is that when using multiscale deep learning method, we utilize the case 100 dimensions to replace that of 20 dimensions. All of the numerical experiments have been performed in Python 3.8 using TensorFlow 2.4, on NVIDIA Tesla P100 GPU (16GB memory). The simulation codes of proposed method are available on the GitHub page \footnote{\url{https://github.com/xiaoxu1996/Deep-PDEs}}.

\subsection{\normalfont High-dimensional AC equation}\label{Se4.1}
\begin{figure}[H]
\centering
\label{F1.1}
\begin{minipage}[b]{0.48\textwidth}
\includegraphics[width=\textwidth,height=0.25\textheight]{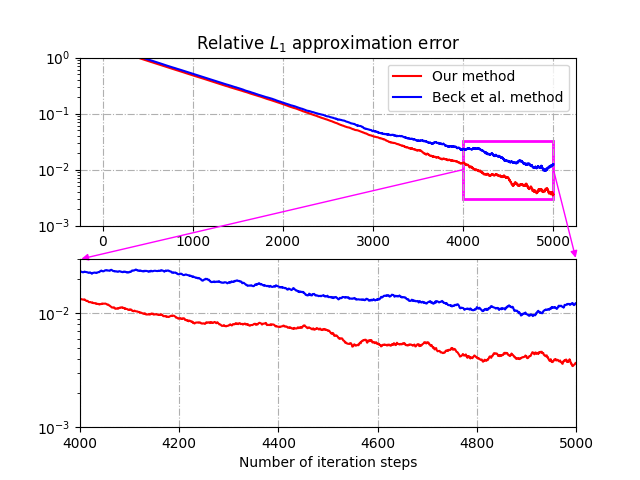}
\end{minipage}
\label{F1.2}
\begin{minipage}[b]{0.48\textwidth}
\includegraphics[width=\textwidth,height=0.25\textheight]{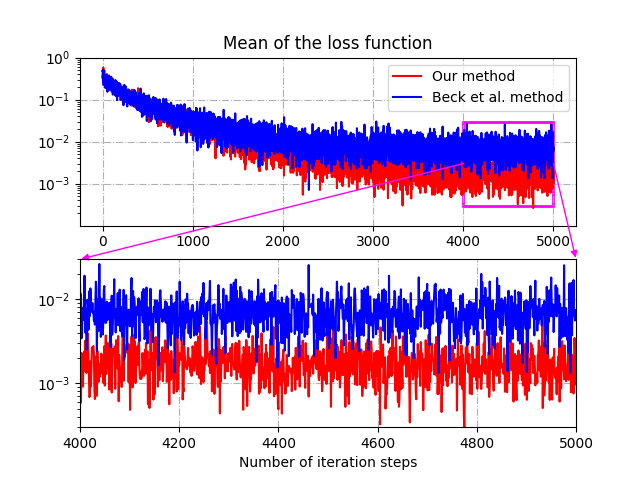}
\end{minipage}
\caption{\quad Relative $L_1$ approximation error and the mean of the empirical loss function of the 20-dimensional AC equation.}
\label{F1}
\end{figure}

This subsection discusses the approximate solution of the high-dimensional AC equation with a cubic nonlinearity (see \eqref{eq4.3}). Next, the following two examples show that the approximated calculation of Allen equations of different dimensions from multiscale deep learning fusion and convolutional neural networks, respectively.
\begin{example}\normalfont \label{ex1}
\textit{\bf Multiscale deep learning fusion}.
Assuming the notations $T=\frac{3}{10}$, $\widetilde{\gamma}=\frac{1}{1000}, d=20,\tilde{d}\in\{20,30,40,50\},N=20,\xi=\{0,\ldots,0\}\in\mathbb{R}^{d},t\in[0,T),\mathbf{x},\mathbf{z}\in\mathbb{R}^d,\mathbf{y}\in\mathbb{R},S\in\mathbb{R}^{d\times d},t_s=\frac{sT}{N}$, $\hat{g}(\mathbf{x})=[2+\frac{2}{5}||\mathbf{x}||_{\mathbb{R}^d}^2]^{-1}$, and
\begin{equation}
\label{eq4.1}
f(t,\mathbf{x},\mathbf{y},\mathbf{z},S)=-\frac{1}{2}{\rm Trace}(S)-\mathbf{y}+\mathbf{y}^3,
\end{equation}
and suppose that $u:[0,T]\times\mathbb{R}^d\to\mathbb{R}$ is an at most polynomially growing continuous function, such that $u(T,\mathbf{x})=\hat{g}(\mathbf{x}),u|_{[0,T)\times\mathbb{R}^d}\in
C^{1,3}([0,T)\times\mathbb{R}^d,\mathbb{R})$, and
\begin{equation}
\label{eq4.2}
\frac{\partial u}{\partial t}(t,\mathbf{x})=f(t,\mathbf{x},u(t,\mathbf{x}),(\nabla_\mathbf{x} u)(t,\mathbf{x}),({\rm Hess}_\mathbf{x} u)(t,\mathbf{x})),
\end{equation}
for all $(t,\mathbf{x})\in[0,T)\times\mathbb{R}^d$. The solution $u:[0,T)\times\mathbb{R}^d\to\mathbb{R}$ of \eqref{eq4.2} such that $u(T,\mathbf{x})=\left[2+\frac{2}{5}||\mathbf{x}||_{\mathbb{R}^d}\right]^{-1}$ and
\begin{equation}
\label{eq4.3}
\frac{\partial u}{\partial t}(t,\mathbf{x})+\frac{1}{2}(\triangle_\mathbf{x} u)(t,\mathbf{x})+u(t,\mathbf{x})-[u(t,\mathbf{x})]^3=0,
\end{equation}
for all $(t,\mathbf{x})\in [0,T)\times\mathbb{R}^d$.
\end{example}
Table \ref{table1} displays different methods to approximatively calculate the mean and standard deviation of $u^{\Theta_m}$ (i.e., $\mu_{u^{\Theta_m}}$ and $\sigma_{u^{\Theta_m}}$), the mean and standard deviation of corresponding $L_1$-approximation error associated to $u^{\Theta_m}$ (i.e., $\mu_{L^1_{\rm error}}$ and $\sigma_{L^1_{\rm error}}$), and the runtime in seconds needed to calculate one realization of $u^{\Theta_m}$ against $m\in\{0,1000,2000,3000,4000,5000\}$ based on 10 independent runs. In addition, Figure \ref{F1} depicts approximations of the mean of the relative $L_1$-approximation error and approximations of the mean of the loss function associated to $u^{\Theta_m}$ against $m\in\{0,1,2,\ldots,5000\}$ based on 10 independent realizations. In the approximative calculations of the relative $L^1$-approximation error, the value $u(0,\xi)$ of the solution $u$ of the \eqref{eq4.3} has been replaced by the value 0.30879 which, in turn, has been calculated through the Branching diffusion method \cite{33}. In particular, the relative $L_1$-approximation error is calculated as $\frac{|u^{\Theta_m}-0.30879|}{0.30879}$.

It is not difficult to see from Table \ref{table1} that the approximate solution obtained by our method has higher accuracy, and the running time is also greatly reduced. To more intuitively compare with the existing methods, we draw Figure \ref{F1}. Regarding the relative $L^1$-approximation error in Figure \ref{F1}, proposed method is almost consistent with the method in \cite{Beck} when the number of training steps is small. However, as the number of training steps increases, the proposed method has a smaller relative $L^1$-approximation error, which means that our method is more accurate and effective. In addition, we purposely magnify the relative $L^1$-approximation error from steps 4000 to 5000 to the lower part of the figure. From the enlarged picture, it can be clearly seen that our relative $L^1$-approximation error is already less than 0.01. At the same time, the right side of Figure \ref{F1} shows the trend of the loss function. As shown, our loss function value is smaller.

\begin{table}[H]
\begin{center} \small
\caption{Numerical simulations of the 20-dimensional AC equation.}
\label{table1}
\begin{tabular}{|c|m{1.3cm}<{\centering}|c|c|c|c|m{1.2cm}<{\centering}|m{1.4cm}<{\centering}|}
\hline
Method&Training steps&$\displaystyle\mu_{u^{\Theta_m}}$&$\displaystyle\sigma_{u^{\Theta_m}}$&$\displaystyle\mu_{L^1_{\rm error}}$&$\displaystyle\sigma_{L^1_{\rm error}}$
&Mean of the loss function&Runtime in sec.\\
\hline
{Beck et al. \cite{Beck}}
&0   &-0.04958&0.57116&1.88360&1.10466&0.47839&6\\
&1000&0.19091 &0.14298&0.51528&0.30760&0.02459&14\\
&2000&0.26892 &0.04361&0.15655&0.11004&0.01089&23\\
&3000&0.29646 &0.01359&0.04874&0.03397&0.00724&31\\
&4000&0.30252 &0.00584&0.02369&0.01444&0.01550&40\\
&5000&0.30584 &0.00288&0.01243&0.00487&0.00662&49\\
\hline
{Our results}
&0   &-0.02988&0.58509&1.78238&1.27133&0.35253&2\\
&1000&0.20342 &0.15110&0.48308&0.35003&0.01850&3\\
&2000&0.27478 &0.04546&0.14750&0.10976&0.00412&5\\
&3000&0.29954 &0.01301&0.03965&0.03319&0.00139&6\\
&4000&0.30582 &0.00393&0.01328&0.00881&0.00120&7\\
&5000&0.30852 &0.00123&0.00363&0.00184&0.00232&9\\
\hline
\end{tabular}
\end{center}
\end{table}

\begin{example}\normalfont\label{ex2}
\textit{\bf Convolutional neural networks}. We still utilize certain basic settings from Example \ref{ex1}, and the only thing that needs to be modified is the dimension of the data. Here, set $d=256$ or $d=400$.
\begin{table}[H]
\begin{center} \small
\caption{Numerical simulations of the large-dimensional AC equation.}
\label{table2}
\begin{tabular}{|c|m{1.3cm}<{\centering}|c|c|c|c|m{1.2cm}<{\centering}|m{1.4cm}<{\centering}|}
\hline
Dimension&Training steps&$\displaystyle\mu_{u^{\Theta_m}}$&$\displaystyle\sigma_{u^{\Theta_m}}$&$\displaystyle\mu_{L^1_{\rm error}}$&$\displaystyle\sigma_{L^1_{\rm error}}$
&Mean of the loss function&Runtime in sec.\\
\hline
{$d=256$}
&0    &-0.15151&0.57393&12.7026&7.15709&0.74728&2\\
&2000 &0.03103 &0.03149&0.73467&0.31448&0.02349&4\\
&4000 &0.04045 &0.00364&0.06574&0.06365&0.00616&7\\
&6000 &0.04217 &0.00131&0.02855&0.02033&0.00087&10\\
&8000 &0.04139 &0.00042&0.00797&0.00723&0.00010&12\\
&10000&0.04155 &0.00011&0.00227&0.00158&0.00003&15\\
\hline
{$d=400$}
&0    &0.08637&0.46341&14.6610&9.05968&0.32806&2\\
&2000 &0.02730&0.01902&0.59862&0.36606&0.04361&4\\
&4000 &0.02499&0.00474&0.16069&0.10427&0.00721&7\\
&6000 &0.02685&0.00161&0.05114&0.03197&0.00239&10\\
&8000 &0.02698&0.00082&0.02202&0.02099&0.00028&13\\
&10000&0.02729&0.00022&0.00850&0.00619&0.00004&15\\
\hline
\end{tabular}
\end{center}
\end{table}

\begin{figure}[H]
\centering
\label{F2.1}
\begin{minipage}[b]{0.48\textwidth}
\includegraphics[width=\textwidth,height=0.25\textheight]{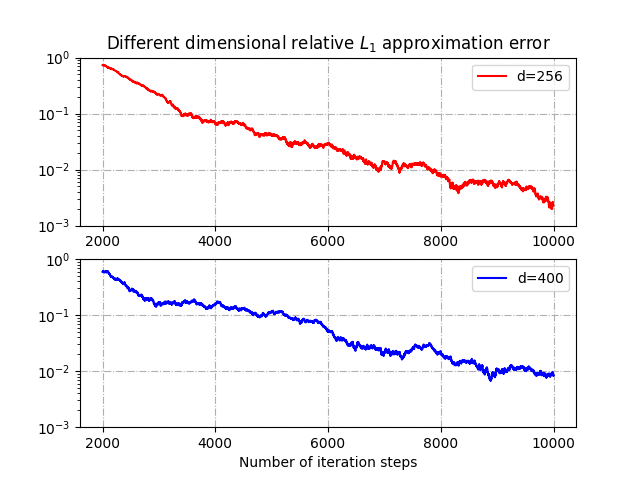}
\end{minipage}
\label{F2.2}
\begin{minipage}[b]{0.48\textwidth}
\includegraphics[width=\textwidth,height=0.25\textheight]{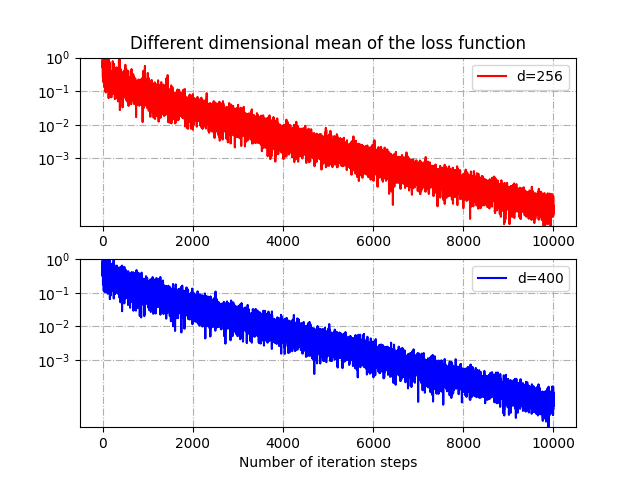}
\end{minipage}
\caption{\quad  Relative $L_1$ approximation error and the mean of the empirical loss function of the large-dimensional AC equation.}
\label{F2}
\end{figure}
\end{example}
 Table \ref{table2} extracts approximate solutions of $u^{\Theta_m}$ in different dimensions by convolutional neural networks. The difference with Example \ref{ex1} is that the number of iteration steps here $m\in\{0,2000,4000,6000,8000,10000\}$. And in Figure \ref{F2}, $m\in\{0,1,2,\ldots,10000\}$. Besides, the 256- and 400-dimensional value $u(0,\xi)$ of the solution $u$ of the \eqref{eq4.3} has been replaced by the value 0.041531 and 0.027106, which is also calculated through the Branching diffusion method \cite{33}. Hence, the different dimensional relative $L_1$-approximation error is calculated as $\dfrac{|u^{\Theta_m}-0.041531|}{0.041531},\dfrac{|u^{\Theta_m}-0.027106|}{0.027106}$, respectively.

In Table \ref{table2} and Figure \ref{F2}, no matter whether the dimension of the equation is 256 or 400, as the number of iteration steps increases, the relative $L_1$-approximation error of the approximate solution decreases gradually, and the loss function also tends to decrease in general. This shows that it is numerically feasible for us to use convolutional neural networks to approximately solve higher-dimensional stochastic PDEs.

\subsection{\normalfont High-dimensional BSB equation}\label{Se4.2}
This subsection presents  the calculation of the high-dimensional BSB equation (see \cite{2} and \eqref{eq4.10}). Similarly, we employ two examples to show that.

\begin{example}\normalfont \label{ex5}
\textit{\bf Multiscale deep learning fusion}.
Suppose that $T=1, d=100,\tilde{d}\in\{75,100,50,125\},N=20,\epsilon=10^{-8}$, and assume for all $\omega\in\Omega$ that $\xi(\omega)=(1,\frac{1}{2},1,\frac{1}{2},\ldots,1,\frac{1}{2})\in\mathbb{R}^d$. Set
\begin{equation}
\label{eq4.7}
\widetilde{\gamma}_m=1.0\cdot\left(\frac{1}{2}\right)^{[m/200]}.
\end{equation}
Here, $[\cdot]$ represents taking the integer of $m/200$. By setting
$\sigma_{max}=\frac{4}{10}$,  $\sigma_{min}=\frac{1}{10}$,  $\sigma_c=\frac{4}{10}$, let us  define  the function $\bar{\sigma}:\mathbb{R}\to\mathbb{R}$ as
\begin{equation}
\label{eq4.8}
\bar{\sigma}(x)=\left\{
\begin{array}{cc}
\sigma_{max},&x\geq0,\\
\sigma_{min},&x<0
\end{array}
\right.
\end{equation}
for all $x\in\mathbb{R}$. Assuming for all $s,t\in[0,T]$, $\mathbf{x}=(\mathbf{x}_1,\ldots,\mathbf{x}_d)$, $\mathbf{w}=(\mathbf{w}_1,\ldots,\mathbf{w}_d)$,
$\mathbf{z}=(\mathbf{z}_1,\ldots,\mathbf{z}_d)\in\mathbb{R}^d,y\in\mathbb{R}$, $S=(S_{ij})_{(i,j)\in\{1,\ldots,d\}^2}\in\mathbb{R}^{d\times d}$, we have that $\sigma(\mathbf{x})=\sigma_c \text{diag}(\mathbf{x}_1,\ldots,\mathbf{x}_d)$, $\mathcal{H}(s,t,\mathbf{x},\mathbf{w})=\mathbf{x}+\sigma(\mathbf{x})\mathbf{w}$, $\hat{g}(\mathbf{x})=||\mathbf{x}||_{\mathbb{R}^d}^2$, and that
\begin{equation}
\label{eq4.9}
f(t,\mathbf{x},\mathbf{y},\mathbf{z},S)=-\frac{1}{2}\sum_{i=1}^d |\mathbf{x}_i|^2|\bar{\sigma}(S_{ii})|^2 S_{ii}+\hat{r}(\mathbf{y}-\langle \mathbf{x},\mathbf{z}\rangle_{\mathbb{R}^d}).
\end{equation}
The solution $u:[0,T]\times\mathbb{R}^d\to\mathbb{R}$ such that $u(T,\mathbf{x})=||\mathbf{x}||_{\mathbb{R}^d}^2$ and
\begin{equation}
\label{eq4.10}
\frac{\partial u}{\partial t}(t,\mathbf{x})+\frac{1}{2}\sum_{i=1}^d|\mathbf{x}_i|^2|\bar{\sigma}(\frac{\partial^2 u}{\partial \mathbf{x}_i^2}(t,\mathbf{x}))|^2\frac{\partial^2 u}{\partial \mathbf{x}_i^2}(t,\mathbf{x})=\hat{r}(u(t,\mathbf{x})-\langle \mathbf{x},(\nabla_\mathbf{x} u)(t,\mathbf{x})\rangle_{\mathbb{R}^d})
\end{equation}
for all $(t,\mathbf{x})\in[0,T)\times\mathbb{R}^d$.
\end{example}

\begin{figure}[H]
\centering
\label{F5.1}
\begin{minipage}[b]{0.48\textwidth}
\includegraphics[width=\textwidth,height=0.25\textheight]{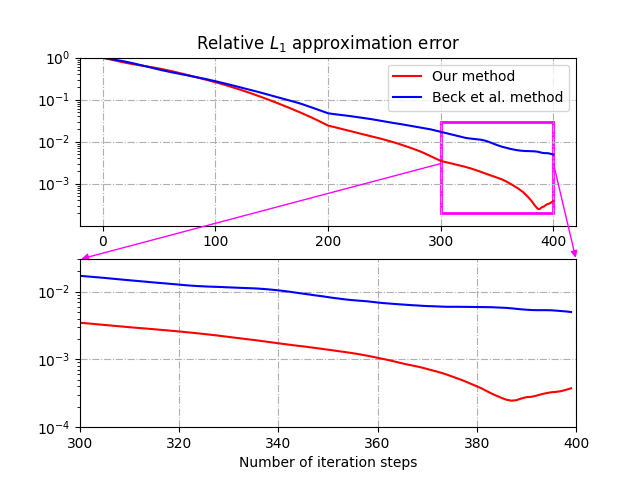}
\end{minipage}
\label{F5.2}
\begin{minipage}[b]{0.48\textwidth}
\includegraphics[width=\textwidth,height=0.25\textheight]{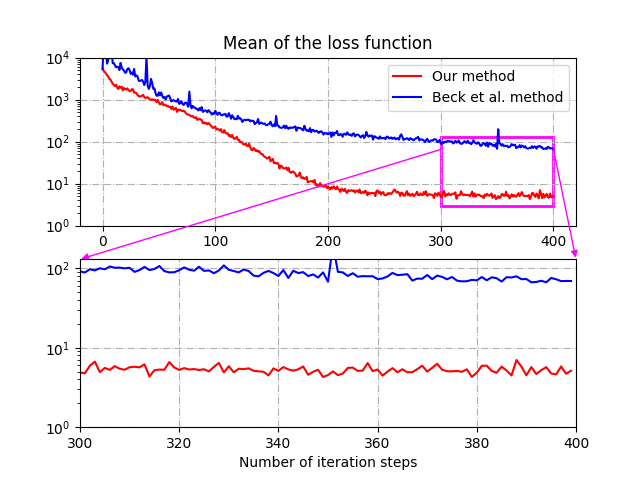}
\end{minipage}
\caption{\quad  Relative $L_1$ approximation error and the mean of the empirical loss function of the 100-dimensional BSB equation.}
\label{F5}
\end{figure}

Table \ref{table5} lists different methods to approximatively calculate the mean and standard deviation of  $u^{\Theta_m}$, the mean and standard deviation of corresponding $L_1$-approximation error associated to $u^{\Theta_m}$, and the runtime in seconds, needed to calculate one realization of $u^{\Theta_m}$ against $m\in\{0,100,200,300,400\}$ based on 10 independent runs. In addition, Figure \ref{F5} depicts approximations of the mean of the relative $L_1$-approximation error and approximations of the mean of the loss function associated to $u^{\Theta_m}$ against $m\in\{0,1,2,\ldots,400\}$ based on 10 independent realizations. In the approximative calculations of the relative $L^1$-approximation error, the value $u(0,(1,\frac{1}{2},1,\frac{1}{2},\ldots,1,\frac{1}{2}))$ of the solution $u$ of \eqref{eq4.10} has been replaced by the value 77.1049, in turn, which has been calculated by means of Lemma \ref{lemma4.1} below (more details see \cite{Beck}). The relative $L_1$-approximation error is $\frac{|u^{\Theta_m}-77.1049|}{77.1049}$.

\begin{lemma}\label{lemma4.1}\normalfont
Suppose that $0<c,\sigma_{max},r,T<\infty$, $0<\sigma_{min}<\sigma_{max}$, $d\in\mathbb{N}$, and assume $\bar{\sigma}:\mathbb{R}\to\mathbb{R}$ is the function, such that
\begin{equation}
\label{eq4.11}
\bar{\sigma}(x)=\left\{
\begin{array}{cc}
\sigma_{max},&x\geq0,\\
\sigma_{min},&x<0
\end{array}
\right.
\end{equation}
for all $x\in\mathbb{R}$, and we let $\hat{g}:\mathbb{R}^d\to\mathbb{R}$ and $u:[0,T]\times\mathbb{R}^d\to\mathbb{R}$ be the functions, such that $\hat{g}(\mathbf{x})=c||\mathbf{x}||_{\mathbb{R}^d}^2=c\sum_{i=1}^d|\mathbf{x}_i|^2$ and
\begin{equation}
\label{eq4.12}
u(t,\mathbf{x})=\exp([r+|\sigma_{max}|^2](T-t))\hat{g}(\mathbf{x})
\end{equation}
for all $t\in[0,T]$, $\mathbf{x}=(\mathbf{x}_1,\ldots,\mathbf{x}_d)\in\mathbb{R}^d$. Then, we have for all $t\in[0,T],\mathbf{x}=(\mathbf{x}_1,\ldots,\mathbf{x}_d)\in\mathbb{R}^d$ that $u\in C^{\infty}([0,T]\times\mathbb{R}^d,\mathbb{R}),u(T,\mathbf{x})=\hat{g}(\mathbf{x})$, and
\begin{equation}
\label{eq4.13}
\frac{\partial u}{\partial t}(t,\mathbf{x})+\frac{1}{2}\sum_{i=1}^d|\mathbf{x}_i|^2|\bar{\sigma}(\frac{\partial^2 u}{\partial \mathbf{x}_i^2}(t,\mathbf{x}))|^2\frac{\partial^2 u}{\partial \mathbf{x}_i^2}(t,\mathbf{x})=\hat{r}(u(t,\mathbf{x})-\langle \mathbf{x},(\nabla_\mathbf{x} u)(t,\mathbf{x})\rangle_{\mathbb{R}^d}).
\end{equation}
\end{lemma}
Looking at Table \ref{table5} as a whole  we observe that the approximate solution obtained by our method has higher accuracy. However, unlike Example \ref{ex1}, our runtime will be a bit more. Similarly, we paint Figure \ref{F5} for comparing with the existing methods. It is evident from Figure \ref{F5} that when the number of iteration steps exceeds 200, the proposed method already stratifies with the method of Beck et al. \cite{Beck}. And from the partially enlarged picture, Beck et al. \cite{Beck} method differs from us by one coordinate scale in terms of the relative $L^1$-approximation error and loss function value. These all demonstrate and illustrate the effectiveness of our method.

\begin{table}[H]
\begin{center} \small
\caption{Numerical simulations of the 100-dimensional BSB equation.}
\label{table5}
\begin{tabular}{|c|m{1.5cm}<{\centering}|c|c|c|c|m{2cm}<{\centering}|m{1.5cm}<{\centering}|}
\hline
Method&Training steps&$\displaystyle\mu_{u^{\Theta_m}}$&$\displaystyle\sigma_{u^{\Theta_m}}$&$\displaystyle\mu_{L^1_{\rm error}}$&$\displaystyle\sigma_{L^1_{\rm error}}$
&Mean of the loss function&Runtime in sec.\\
\hline
{Beck et al. \cite{Beck}}
&0   &0.3940&0.2253&0.99489&0.00292&5355.51&23\\
&100&55.9301&1.9195&0.27462&0.02489&540.55 &27\\
&200&73.4561&0.9547&0.04732&0.01238&149.26 &31\\
&300&75.7877&0.5027&0.01708&0.00652&90.979 &36\\
&400&76.7701&0.3009&0.00491&0.00316&63.846 &40\\
\hline
{Our results}
&0   &0.5517&0.2378&0.99285&0.00308&5411.35&21\\
&100&57.0542&0.4246&0.26004&0.00551&226.89 &29\\
&200&75.2420&0.1476&0.02416&0.00191&8.619  &36\\
&300&76.8373&0.0515&0.00347&0.00067&4.866  &44\\
&400&77.1226&0.0302&0.00039&0.00024&4.882  &52\\
\hline
\end{tabular}
\end{center}
\end{table}

\begin{figure}[H]
\centering
\label{F6.1}
\begin{minipage}[b]{0.48\textwidth}
\includegraphics[width=\textwidth,height=0.25\textheight]{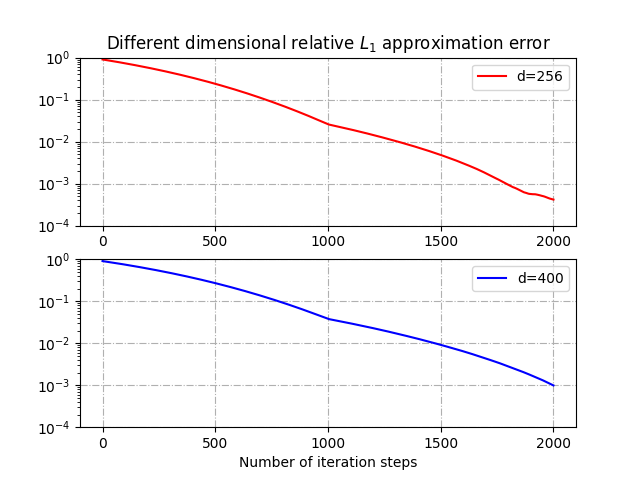}
\end{minipage}
\label{F6.2}
\begin{minipage}[b]{0.48\textwidth}
\includegraphics[width=\textwidth,height=0.25\textheight]{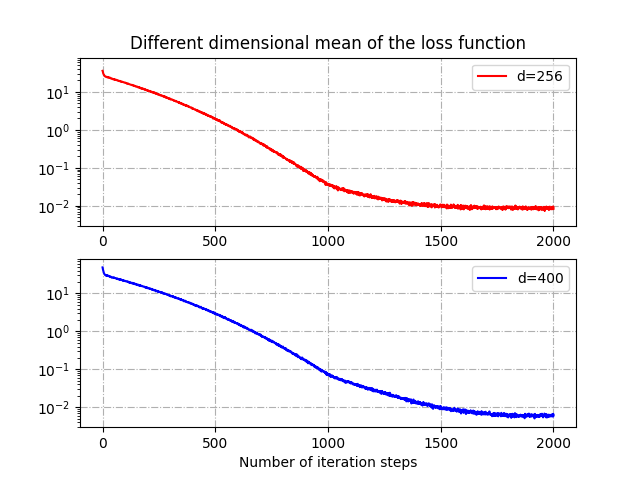}
\end{minipage}
\caption{\quad  Relative $L_1$ approximation error and the mean of the empirical loss function of the large-dimensional BSB equation.}
\label{F6}
\end{figure}

\begin{example}\normalfont \label{ex6}
\textit{\bf Convolutional neural networks}. Herein, most of our settings are the same as Example \ref{ex5}. Based on this point, what needs to be modified is the dimension of the data and the learning rate. Firstly, we set $d=256$ or $d=400$, and the learning rate is
\begin{equation}
\label{eq4.14}
\widetilde{\gamma}_m=2.0\cdot\left(\frac{1}{2}\right)^{[m/500]}.
\end{equation}
\end{example}
Table \ref{table6} reports approximate solutions of $u^{\Theta_m}$ in different dimensions by convolutional neural networks. The difference with Example \ref{ex5} is that the number of iteration steps here $m\in\{0,200,400,600,800,1000\}$. In addition, in Figure \ref{F6}, $m\in\{0,1,2,\ldots,1000\}$. Also, the 256- and 400-dimension value $u\left(0,(1,\frac{1}{2},1,\frac{1}{2},\ldots,1,\frac{1}{2})\right)$ of the solution $u$ of \eqref{eq4.10} has been replaced via the value 197.3885 and 308.4195, respectively. It also can be computed by means of Lemma \ref{lemma4.1}. And the different dimensions relative $L_1$-approximation error is $\frac{|u^{\Theta_m}-197.3885|}{197.3885}$, $\frac{|u^{\Theta_m}-308.4195|}{308.4195}$, respectively.

In Table \ref{table6}, it can be seen that from $256$ dimensions to $400$ dimensions, the running time using convolutional neural networks increases exponentially. This is mainly because as the dimension increases, the memory overhead increases. However, the accuracy of the approximated solution did not change much. This demonstrates that convolutional neural networks can extend approximated solutions to higher dimensions without losing accuracy. Also, Figure \ref{F6} can show this more intuitively.

\begin{table}[H]
\begin{center} \small
\caption{Numerical simulations of the large-dimensional BSB equation.}
\label{table6}
\begin{tabular}{|c|m{1.5cm}<{\centering}|c|c|c|c|m{2cm}<{\centering}|m{1.5cm}<{\centering}|}
\hline
Dimension&Training steps&$\displaystyle\mu_{u^{\Theta_m}}$&$\displaystyle\sigma_{u^{\Theta_m}}$&$\displaystyle\mu_{L^1_{\rm error}}$&$\displaystyle\sigma_{L^1_{\rm error}}$
&Mean of the loss function&Runtime in sec.\\
\hline
{$d=256$}
&0   &0.4901  &0.2948&0.99752&0.00149&35095 &4\\
&200 &164.3867&0.6826&0.16719&0.00346&345.59&34\\
&400 &190.3597&0.2866&0.03561&0.00145&26.737&63\\
&600 &194.5438&0.1614&0.01441&0.00082&17.643&92\\
&800 &196.7375&0.1106&0.00330&0.00056&15.022&122\\
&1000&197.3413&0.0793&0.00041&0.00023&14.395&151\\
\hline
{$d=400$}
&0   &0.5218  &0.2603&0.99831&0.00084&86234 &6\\
&200 &170.3224&1.6187&0.44776&0.00525&3696.6&77\\
&400 &271.3132&1.1839&0.12031&0.00384&347.57&148\\
&600 &298.8917&0.6204&0.03089&0.00201&42.596&219\\
&800 &305.8975&0.2641&0.00818&0.00086&29.087&291\\
&1000&308.5768&0.1068&0.00051&0.00035&23.190&362\\
\hline
\end{tabular}
\end{center}
\end{table}

\subsection{\normalfont High-dimensional HJB equation}\label{Se4.3}
\begin{figure}[H]
\centering
\label{F3.1}
\begin{minipage}[b]{0.48\textwidth}
\includegraphics[width=\textwidth,height=0.25\textheight]{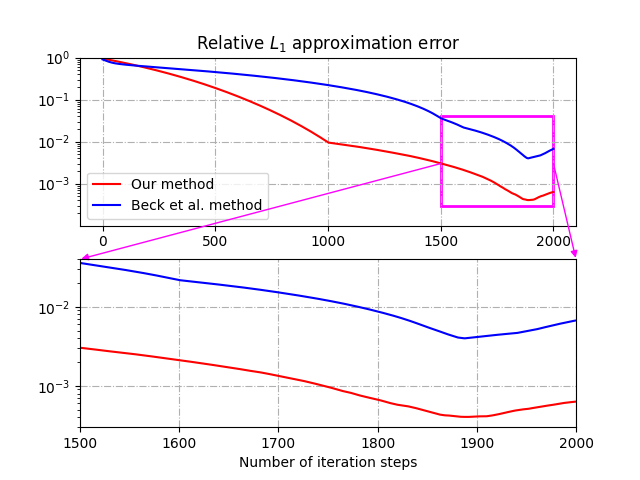}
\end{minipage}
\label{F3.2}
\begin{minipage}[b]{0.48\textwidth}
\includegraphics[width=\textwidth,height=0.25\textheight]{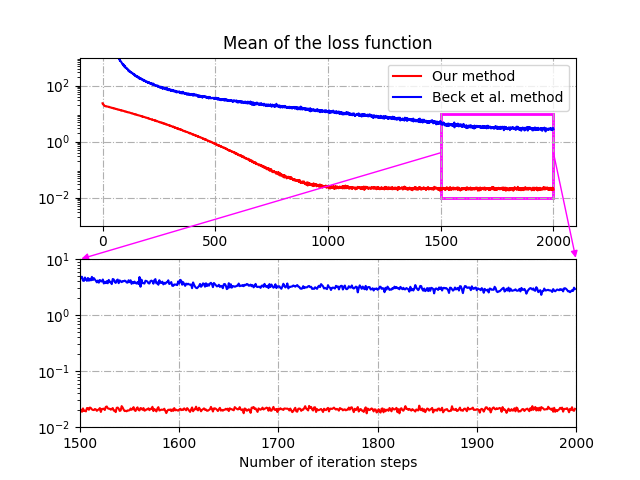}
\end{minipage}
\caption{\quad  Relative $L_1$ approximation error and the mean of the empirical loss function of the 100-dimensional HJB equation.}
\label{F3}
\end{figure}

This subsection approximatively calculates the solution of a high-dimensional HJB equation with a nonlinearity that is quadratic in the gradient (see \cite{33}). In the following,  we present two examples to show the related calculation.
\begin{example}\normalfont \label{ex3}
\textit{\bf Multiscale deep learning fusion}.
We suppose $d=100,\tilde{d}\in\{50,75,100,125\},T=1,N=20,\epsilon=10^{-8}$, and suppose for all $\omega\in\Omega$ that
$\xi(\omega)=\boldsymbol{0}\in\mathbb{R}^d$. Then assume for all $m\in\mathbb{N}_0,s,t\in[0,T],\mathbf{x},\mathbf{w},\mathbf{z}\in\mathbb{R}^d,\mathbf{y}\in\mathbb{R},S\in\mathbb{R}^{d\times d}$ that $\sigma(\mathbf{x})=\sqrt{2}\mathrm{Id}_{\mathbb{R}^d},\mathcal{H}(s,t,\mathbf{x},\mathbf{w})=\mathbf{x}+\sqrt{2}\mathbf{w},\hat{g}(\mathbf{x})
=\ln\left(\frac{1}{2}[1+||\mathbf{x}||_{\mathbb{R}^d}^2]\right),$ $f(t,\mathbf{x},\mathbf{y},\mathbf{z},S)=-{\rm Trace}(S)-||\mathbf{z}||_{\mathbb{R}^d}^2$, and
\begin{equation}
\label{eq4.4}
\widetilde{\gamma}_m=\frac{1}{100}\cdot\left(\frac{1}{5}\right)^{[m/1000]}.
\end{equation}
The solution $u:[0,T)\times\mathbb{R}^d\to\mathbb{R}$ of the PDE \eqref{eq4.2} satisfies for all $(t,\mathbf{x})\in
[0,T)\times\mathbb{R}^d$ that
\begin{equation}
\label{eq4.5}
\frac{\partial u}{\partial t}(t,\mathbf{x})+(\triangle_\mathbf{x} u)=||\nabla_\mathbf{x} u(t,\mathbf{x})||_{\mathbb{R}^d}^2.
\end{equation}
\end{example}
Table \ref{table3} lists different methods to approximatively calculate the mean and standard deviation of $u^{\Theta_m}$, the mean and standard deviation of relative $L_1$-approximation error associated to $u^{\Theta_m}$, and the runtime in seconds, needed to calculate one realization of $u^{\Theta_m}$ against $m\in\{0,500,1000,1500,2000\}$, based on 10 independent runs. Furthermore, Figure \ref{F3} shows approximations of the mean of the relative $L_1$-approximation error and approximations of the mean of the loss function associated to $u^{\Theta_m}$ against $m\in\{0,1,2, \ldots, 2000\}$ based on 10 independent realizations. For the approximative calculations of the relative $L_1$-approximation error, the value $u(0,\xi)$ of the solution $u$ of \eqref{eq4.5} has been substituted by the value 4.5901, conversely, which was calculated by the means of in \cite[Lemma 4.2]{33} and the classical Monte Carlo method \cite{33}.

It can be clearly observed from Table \ref{table3} and Figure \ref{F3} that the approximated solution obtained via our method has higher accuracy.
 Figure \ref{F3}, the curve slope of the relative $L_1$-approximation error and the loss function change with our method at 1000 steps, which is mainly caused by the change of the learning rate (see \eqref{eq4.4}). Likewise, we also place the local comparison from steps 1500 to 2000 at the bottom of this figure.

\begin{table}[H]
\begin{center} \small
\caption{Numerical simulations of the 100-dimensional HJB equation.}
\label{table3}
\begin{tabular}{|c|m{1.5cm}<{\centering}|c|c|c|c|m{2cm}<{\centering}|m{1.5cm}<{\centering}|}
\hline
Method&Training steps&$\displaystyle\mu_{u^{\Theta_m}}$&$\displaystyle\sigma_{u^{\Theta_m}}$&$\displaystyle\mu_{L^1_{\rm error}}$&$\displaystyle\sigma_{L^1_{\rm error}}$
&Mean of the loss function&Runtime in sec.\\
\hline
{Beck et al. \cite{Beck}}
&0   &0.4328&0.0620&0.90571&0.01351&1065.5&17\\
&500 &2.5108&0.0555&0.45300&0.01208&37.574&33\\
&1000&3.5726&0.0432&0.22168&0.00942&11.839&49\\
&1500&4.4255&0.0293&0.03587&0.00639&5.105 &65\\
&2000&4.6101&0.0258&0.00673&0.00232&2.783 &81\\
\hline
{Our results}
&0   &0.2294&0.0940&0.95001&0.02047&23.32&18\\
&500 &3.7223&0.0603&0.18907&0.01313&0.834&42\\
&1000&4.5465&0.0097&0.00951&0.00212&0.025&67\\
&1500&4.5762&0.0052&0.00304&0.00113&0.022&91\\
&2000&4.5924&0.0021&0.00063&0.00024&0.019&115\\
\hline
\end{tabular}
\end{center}
\end{table}

\begin{figure}[H]
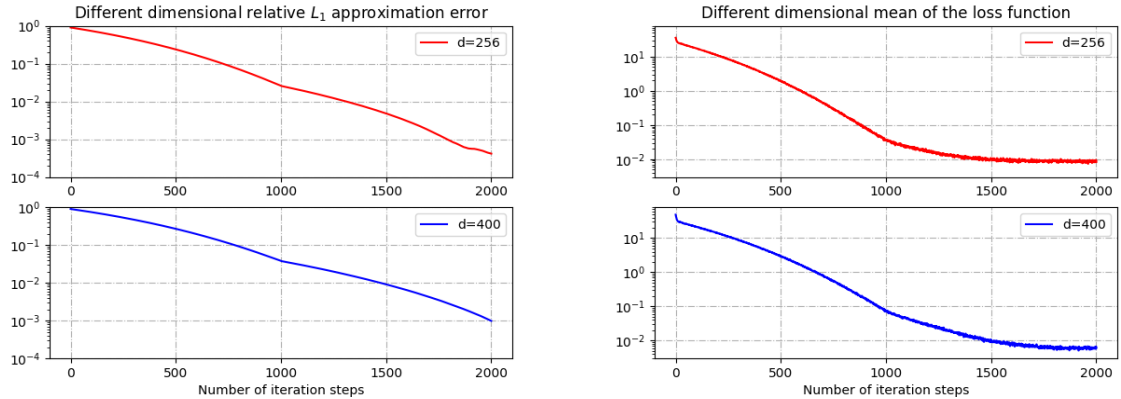

\centering
\label{F4.1}
\begin{minipage}[b]{0.48\textwidth}
\includegraphics[width=\textwidth,height=0.25\textheight]{H_3.png}
\end{minipage}
\label{F4.2}
\begin{minipage}[b]{0.48\textwidth}
\includegraphics[width=\textwidth,height=0.25\textheight]{H_4.png}
\end{minipage}
\caption{\quad  Relative $L_1$ approximation error and the mean of the empirical loss function of the large-dimensional HJB equation.}
\label{F4}
\end{figure}

\begin{example}\normalfont \label{ex4}
\textit{\bf Convolutional neural networks}. Herein, certain basic settings from Example \ref{ex3} are still used, and the only thing that needs to be changed is the dimension of the data. Below, set $d=256$ or $d=400$. Noting the learning rate, we adjusted the learning rate with a fixed number of steps instead of exponential decay. The specific formula is
\begin{equation}
\label{eq4.6}
\widetilde{\gamma}_m=\left\{
\begin{array}{cc}
0.01,&m<1000,\\
0.005,&m\geq1000.
\end{array}
\right.
\end{equation}
\end{example}
Table \ref{table4} and Figure \ref{F4} display approximated solutions of $u^{\Theta_m}$ in different dimensions by convolutional neural networks. Besides, the 256- and 400-dimension value $u(0,\xi)$ of the solution $u$ of \eqref{eq4.1} has been replaced by the value 5.5393 and 5.9877, which also can be calculated through the classical Monte Carlo method \cite{33}. Thus, the different dimensions relative $L_1$-approximation error is $\frac{|u^{\Theta_m}-5.5393|}{5.5393},\frac{|u^{\Theta_m}-5.9877|}{5.9877}$, respectively.

Comparing Table \ref{table3} and Table \ref{table4}, one can find the fact that the running time of using convolutional neural network is faster than using linear neural network. Generally speaking, higher-dimensional problems require more memory and take longer to compute. While  in Table \ref{table4}, it only takes 7 seconds to calculate the 400-dimensional HJB equation. In addition, from the relative $L_1$-approximation error and loss function in Figure \ref{F4}, the accuracy of the convolutional neural network is almost the same as that of the linear neural network. This shows that convolutional neural networks are more suitable than linear neural networks for the HJB equation.

\begin{table}[H]
\begin{center} \small
\caption{Numerical simulations of the large-dimensional HJB equation.}
\label{table4}
\begin{tabular}{|c|m{1.5cm}<{\centering}|c|c|c|c|m{2cm}<{\centering}|m{1.5cm}<{\centering}|}
\hline
Dimension&Training steps&$\displaystyle\mu_{u^{\Theta_m}}$&$\displaystyle\sigma_{u^{\Theta_m}}$&$\displaystyle\mu_{L^1_{\rm error}}$&$\displaystyle\sigma_{L^1_{\rm error}}$
&Mean of the loss function&Runtime in sec.\\
\hline
{$d=256$}
&0   &0.5348&0.2753&0.90346&0.04970&35.86&1\\
&500 &4.2221&0.1976&0.23779&0.03567&1.974&3\\
&1000&5.3966&0.0502&0.02576&0.00907&0.037&4\\
&1500&5.5126&0.0148&0.00481&0.00267&0.010&6\\
&2000&5.5399&0.0025&0.00042&0.00019&0.008&7\\
\hline
{$d=400$}
&0   &0.5902&0.2538&0.90143&0.04239&48.68&2\\
&500 &4.3749&0.1916&0.26935&0.03200&2.853&3\\
&1000&5.7611&0.0596&0.03784&0.00995&0.076&4\\
&1500&5.9330&0.0207&0.00913&0.00346&0.010&6\\
&2000&5.9818&0.0042&0.00099&0.00071&0.006&7\\
\hline
\end{tabular}
\end{center}
\end{table}

\section{Summary}\label{Se5}
\noindent This paper developed numerical approximation for high-dimensional fully nonlinear merged PDEs and 2BSDEs based on the deep CNN technique. First, the forward discretization was employed in the time direction, and then  two approximation approaches were adopted in the space direction  by the multi-scale deep learning fusion and the convolutional neural networks, from which, the former is more accurate and efficient than the method of Beck et al. \cite{Beck}; the latter can use matrix arrangement to calculate higher-dimensional fully nonlinear PDEs, such as $d=400$. These were reflected in the numerical experiments. Unfortunately, despite the computational improvement, we are temporarily unable to obtain theoretical results of the proposed methods, which will be further considered by us in the future.
 Following the results a future study will try to apply  a temporal second-order approximation combined with a regularized convolutional neural network \cite{Zeiler} for solving high-dimensional fully nonlinear merged PDEs-2BSDEs system, based on the stochastic pooling.


\section*{Declaration of competing interest}
\noindent The authors have not disclosed any competing interests.

\section*{Data availability}

\noindent No data was used for the research described in the article.



\end{document}